\documentclass[12pt]{amsart}
\usepackage{amsfonts}
\usepackage[dvips]{graphicx}
\usepackage{amsmath,amsthm,amscd}
\usepackage{amssymb,amsmath, amsthm, amsxtra}
\usepackage{times} 

\newtheorem{theo+}    {Theorem}      [section]
\newtheorem{prop+}  [theo+]  {Proposition}
\newtheorem{coro+}  [theo+]  {Corollary}
\newtheorem{lemm+}  [theo+]  {Lemma}
\newtheorem{deep+}  [theo+]  {Deep Result}
\newtheorem{fact+}  [theo+]  {Fact}
\theoremstyle{definition}
\newtheorem{exam+}  [theo+]  {Example}
\newtheorem{rema+}  [theo+]  {Remark}
\newtheorem{defi+}  [theo+]  {Definition}
\newtheorem{xca+}[theo+]{Exercise}

\numberwithin{equation}{section}

\def\beqn{\begin{equation}}
\def\eeqn{\end{equation}}

\def\epf{\qed \enddemo}

\def\del{\delta}
\def\Del{\Delta}

\def\fa{\mathfrak a}

\def\fg{\mathfrak g}
\def\fk{\mathfrak k}
\def\fh{\mathfrak h}
\def\fl{\mathfrak l}

\def\fp{\mathfrak p}
\def\fq{\mathfrak q}

\def\ft{\mathfrak t}

\def\a{\alpha}
\def\b{\beta}

\def\Claminv2{|C(\Lambda)|^{-2}}
\def\ga{\gamma}

\def\varepsi{\varepsilon}
\def\lam{\lambda}

\def\blam{\underline{\bold \lambda}}

\def\ome{\omega}
\def\Ome{\Omega}

\def\Aa2D{A^{\a,2}(D)}
\def\bAa2D{\overline{A^{\a,2}(D)}}
\def\Ab2D{A^{\beta,2}(D)}
\def\bAb2D{\overline{A^{\beta,2}(D)}}

\def\Norm#1_#2{\Vert#1\Vert_{#2}}

\def\2pd#1#2{\frac{\partial^2 #1}{\partial #2^2}}
\def\p11d#1#2#3{\frac{\partial^2 #1}{  \partial #2\partial #3  }}

\def\b{\beta}
\def\ga{\gamma}
\def\Claminv2{|C(\Lambda)|^{-2}}
\def\del{\delta}
\def\Del{\Delta}

\def\varepsi{\varepsilon}
\def\sig{\sigma}
\def\Sig{\Sigma}

\def\lam{\lambda}

\def\ad{\operatorname{ad}}

\def\tanh{\operatorname{tanh}}

\def\exp{\operatorname{exp}}

\def\tanh{\operatorname{tanh}}

\def\Res{\operatorname{Res}}

\def\Aa2D{A^{\a,2}(D)}
\def\bAa2D{\overline{A^{\a,2}(D)}}
\def\Ab2D{A^{\beta,2}(D)}
\def\bAb2D{\overline{A^{\beta,2}(D)}}

\def\m{\underline{\mathbf m}}
\def\n{\underline{\mathbf n}}

\def\ub1#1{\underline{\mathbf 1^{#1}}}

\def\zae{\"a{}}

\def\bc{\mathbb C}
\def\br{\mathbb R}
\def\nat0{\mathbb Z_{\ge 0}} 

\def\bpf{\begin{proof}}
\def\epf{\end{proof}}
\def\beq{\begin{equation}}
\def\eeq{\end{equation}}

\def\bc{\mathbb C}
\def\br{\mathbb R}

\def\draft{\centerline{(Draft {\the \day}/{\the\month} \the \year.)}}

\begin{document}

\def\phila{\phi_{\blam}}
\def\vepm{\varepsilon_{\m, \nu}(\blam)}
\def\vepmp{\varepsilon_{\m^\prime, \nu}(\blam)}
\def\kam{\kappa_{\m, \nu}(\blam)}
\def\nutnu{\pi({\nu})\otimes \overline{\pi(\nu)}}
\def\pia2ta2{\pi({\frac a2})\otimes \overline{\pi(\frac a2)}}
\def\Hnu{\mathcal H_{\nu}}
\def\Fnu{\mathcal F_{\nu}}
\def\Pm{\mathcal P_{\m}}
\def\PV{\mathcal P(V_{\bc})}
\def\Pa{\mathcal P(\mathfrak a)}
\def\Pn{\mathcal P_{\n}}
\def\cL#1nu{\mathcal L_{{#1}, \nu}}
\def\cL#1a2{\mathcal L_{{#1}, \frac a2}}
\def\E#1nu{E_{{#1}, \nu}}
\def\E#1a2{E_{{#1}, \frac a2}}

\def\Dc{D_{\mathbb C}}
\def\Vc{V_{\mathbb C}}
\def\Bnu{B_{\nu}} 
\def\el{e_{\blam}} 
\def\bnul{b_{\nu}(\lam)} 
\def\GGa{\Gamma_{\Ome}} 
\def\fc#1#2{\frac{#1}{#2}} 
\def\SS{\mathcal S}

\title[Branching coefficients]
{Branching coefficients of holomorphic representations
and  Segal-Bargmann transform}

\author{Genkai Zhang}
\address{Department of Mathematics, Chalmers University of Technology
and G\"o{}teborg University,
S-412 96 G\"o{}teborg, Sweden}
\email{genkai@math.chalmers.se}
\thanks{Research supported by the Swedish Natural Science Research Council (NFR)}



\keywords{Holomorphic discrete series,
highest weight  representations, branching rule,
bounded symmetric
domains, real bounded symmetric domains, Jordan pairs,
Jack symmetric polynomials, orthogonal polynomials, Cayley 
identity}

\begin{abstract}
Let $\mathbb D=G/K$ be a 
 complex  bounded symmetric
domain of tube type in a Jordan algebra
$V_{\mathbb C}$, and let $D=H/L =\mathbb D\cap V$
be its real form in a Jordan
algebra $V\subset V_{\mathbb C}$. 
The analytic continuation of the holomorphic
discrete series on $\mathbb D$
forms a family of interesting representations
of $G$. We consider
the restriction on $D$ of the scalar holomorphic
representations of $G$, as a representation of $H$.
The unitary part of the restriction map
gives then a generalization of
the Segal-Bargmann transform. The group
 $L$ is a spherical subgroup of
$K$ and we find a canonical basis
of $L$-invariant polynomials in components
of the Schmid decomposition
and we express them in terms of
the Jack symmetric polynomials.
We prove that the Segal-Bargmann transform of those
$L$-invariant polynomials
are, under the spherical transform on $D$,
multi-variable Wilson type polynomials and we give a simple alternative proof
of their orthogonality relation.
We find  the expansion of the spherical functions
on $D$, when extended to a neighborhood
in $\mathbb D$, in  terms of the $L$-spherical holomorphic polynomials
on $\mathbb D$,
 the coefficients
being the Wilson polynomials. 
\end{abstract}

\maketitle
\baselineskip 1.40pc

\section{Introduction}

Let $\mathbb D=G/K$
be a complex bounded symmetric domain
of tube type.
The weighted Bergman spaces $\mathcal H_{\nu}$ on $\mathbb D$
form unitary representations of $G$
and are also called the scalar holomorphic
discrete series. They have analytic
continuation in terms of the weight $\nu$
and constitute some interesting
and important family of unitary
representations of $G$. It has
turned out that it is very fruitful
 to study the restriction of the holomorphic
representation to certain subgroups,
both from the point of representation
theory and harmonic analysis. In this paper
we will pursue this by studying the
restriction on some real forms of $\mathbb D$.

The domain $\mathbb D$ can be realized
as a unit ball in a Jordan triple $V_{\bc}$.
Let $V$ be a real form of $V_{\bc}$, $V_{\bc}=V+iV$. The
real form $D=V\cap \mathbb D$
is called a real bounded symmetric domain
if the complex conjugation $\tau$
with respect to $V$ keeps $\mathbb D$ invariant.
In this case $D=H/L$
is also a Riemannian symmetric space where
$H$ is a symmetric subgroup of $G$
and $L$ is a symmetric subgroup of $K$.
Thus we have the following
commutative diagram of subgroup inclusions
$$
\begin{CD} 
G @<<<$K$\\
@AAA @AAA\\
$H$ @<<< $L$
\end{CD}
$$

We consider the  branching law of
the holomorphic representation $\Hnu$ on $\mathbb D$
along the diagram. 
The branching of $\Hnu$ under $K$ is given
by the Schmid decomposition, whereas
its restriction to $H$ (the left vertical
line) is given by the generalized Segal-Bargmann
transform (see \cite{oo-restric-1995} and \cite{gkz-bere-rbsd}),
which gives the unitary equivalence
between $\Hnu$ as the $L^2$-space
on the $D$. To continue
this study, we consider the $L$-invariant
elements in $\Hnu$. The branching
along the  lower horizontal line
is then given by the Helgason spherical
transform; so to get diagram around
we need to find all the $L$-invariant
elements in the Schmid decomposition
and calculate their Segal-Bargmann and
spherical transform, this will be
the main task of the present paper.

 Let  $r$ be
the rank of $D$, so that the rank of $\mathbb D$
is also $r$ if the root
system of $H/L$ is of type $A_{r-1}$ or $D_r$,
it is $2r$ other wise; see Appendix 2.
Let $W$ be the Weyl group of the
root system of the symmetric space $H/L$. 
The space  $\mathcal P=\mathcal P(V_{\mathbb C})$
 of holomorphic polynomials on $V_{\mathbb C}$ under the natural
action of the compact group $K$  is decomposed
into irreducible subspaces $\mathcal P_{\n}$
with signature $\n=n_1\ga_1 +\dots +n_r\ga_r$
for types $A$, $D$ or $B$, and
 $\n=n_1\ga_1 +n_1^\prime\ga_1^\prime+
\dots +n_r\ga_r +n_r^\prime\ga_r^\prime$
for other types, with multiplicity one. This decomposition
can also be viewed as the diagonalization
of a system of Cayley-Capelli type operators,
namely $\Del(z)^{-\a}\Del(\partial)\Del(z)^{\a +1}$,
where $\Del$ is the determinant polynomial
of the Jordan triple $V_{\mathbb C}$ and $\a$
are nonnegative integers; moreover the eigenvalues
of those operators separate the spaces $\mathcal P_{\n}$.
 We consider its subspace
$\mathcal P_{\n}^L$
of $L$-invariant elements. By using
the Cartan-Helgason theorem
we find those signatures $\n$, to be called
spherical signature,
for which
$\mathcal P_{\n}^L\ne 0$. Our first goal will
be to find those polynomials.

For type A
those polynomials are the well-known
Jack symmetric polynomials and their
norm has been calculated by
Faraut-Koranyi \cite{FK}; there are also
 studied  intensively by combinatorial method
\cite{Macd-book}. We shall study in \cite{doz}
their Segal-Bargmann  transform both in bounded and unbounded
realization in relation to the Laplace-transform. However
for other types those polynomials have not been
determined in representation theory. We will use
the Chevalley restriction theorem and the Dunkl-Cherednik operators
to find them.

Let $\fh=\fq +\fl$ be the Cartan decomposition
of the Lie algebra $\fh$ of $H$ and let
$\fa$ be a maximal abelian subspace of $\fq$.
Denote $\Sig$ the root system of $(\fh, \fa)$
and $W$ the corresponding Weyl group.
The space $\fq$
can be identified with the Jordan triple
$V$. By a well-known
theorem of Chevalley
the restriction map from
$\mathcal P^L$ of $L$-invariant polynomials on $V_{\bc}$
to $\fa$ is an isomorphism
onto the subspace $\mathcal P(\fa)^W$
of $W$-invariants polynomials
in $\mathcal P(\fa)$. Thus for each
$L$-spherical signature $\n$ there exists up
to constants a unique polynomials in 
$\mathcal P_{\n}$. Let $D_\xi$, $\xi\in \fa$,
be the Dunkl operators. We can construct
a family $U_i$ of
commuting operators acting on polynomials on $\fa$.
For root system of type A, C those operators
have been previously constructed by Dunkl \cite{Dunkl-cmp-98}.
We prove that when
acting on $L$-invariant polynomials $p$
and restricted to $\fa$,
$$
\Del(z)^{-\a}\Del(\partial)\Del(z)^{\a +1}p(z)
=\prod_{i=1}^r (U_i+\a)p(z).
$$
Thus the problem of diagonalizing the Cayley-Capelli operator
is reduced to that for the operators $U_i$. 
 We  find  those
polynomials $p_{\n}\in \mathcal P_{\n}^L$ in terms of the Jack polynomials
and we calculated the Fock space norm of 
those polynomials; for type  C those
are  done in \cite{Dunkl-cmp-98}
by some different method. We can then
find the norm of those polynomials in  the Bergman space 
$\mathcal H_{\nu}$, by using the result of Faraut-Koranyi. 

We calculate then the Fourier and spherical transforms
of the Segal-Bargmann transforms of those
polynomials in the setting of Fock spaces
and respectively Bergman spaces. They are, up
to a factor of the square root of
the symbol of the Berezin transform,
Weyl group invariant orthogonal polynomials
and will be called the \textit{branching or
coupling coefficients}
as  appeared in the title.  In the former case we prove
that they are of the form
$e^{-\frac 1{4\nu} \Vert \blam\Vert^2 }\zeta_{\n, \nu}(\blam)
$, where
$\zeta_{\n, \nu}(\lam)$ are Hermite type polynomials. Let
$J_{\nu}(x, \lam)$ be the Bessel function associated
with the action of $L$ on $V$,
 we prove that in the
expansion of  $e^{\frac{\nu}2 \Vert x\Vert^2 }J_{\nu}(x, \blam)$
in terms of $p_{\n}(x)$
the coefficients are exactly
$\zeta_{\n, \nu}(\blam)$.
In the later case (curved case) the spherical
transform of the Segal-Bargmann transform
of $p_{\n}$ is of the form
$b_{\nu}(\blam)^{\frac 12} \xi_{\n}(\blam)$, where
$b_{\nu}(\blam)$ is the symbol of the Berezin transform;
the symbol  has been found independently
in \cite{Dijk-Pevzner},  \cite{Neretin-beta-int}  and
  \cite{gkz-bere-rbsd}. (See also \cite{UU} for the
case of complex bounded symmetric domains.)
Thus the polynomials $\xi_{\m}(\blam)$ are
$W$-invariant orthogonal polynomial
with respect to the measure $b_{\nu}(\blam)|c(\blam)|^{-2}$
where $c(\blam)$ is the Harish-Chandra $c$-function.
They are some limiting cases of the multi-variable
Wilson 
polynomials  
\cite{vanDiejen-tams}, so that
the measure in the orthogonality
relation
has now an analytic significance. We find
also an expansion of the product
of the spherical function on $D$
with the reproducing kernel in terms
of the $L$-invariant polynomials, the
coefficient being the spherical transform
of their Segal-Bargmann transform (see Theorem 9.1). We give
thus a unification of the
two types of orthogonal polynomials
associated to a root system,
homogeneous symmetric polynomials
of Jack type on one hand and
the non-homogeneous Wilson polynomials on the 
other.
 For a general root system
the Wilson polynomials have been studied
by van-Diejen \cite{vanDiejen-tams}
under certain self-dual condition.

Along our way of
the study we find an isometric version
of the Chevalley restriction theorem
(see Proposition 5.2) and an
analogy of the Capelli identity
expressing the product $\prod x_j$
and  $\prod D_j$ in terms of the Cherednik
operators  $\prod U_j$, which we believe
are also of independent interest;
see \cite{Torossian-res} and
\cite{HU} for the related study.

We remark that  Theorem \ref{norm-pn-c} 
 can be deduced from \cite{Dunkl-cmp-98},
provided that one proves that they
are eigenpolynomials of the (product of) the Cherednik operators
and identify ours with that of Dunkl;
 in that
paper Dunkl studied the a more general class of polynomials
invariant under certain subgroups of the Weyl group
and found the norm.

We refer also the reader to \cite{Howe-spgp} 
for some algebraic consideration about finding
polynomial invariants of linear groups and
 \cite{toshi-multfree-rest}, \cite{toshi-saga}
for some general results concerning
the branching of unitary highest weight representations.
We mention also that our results can also be interpreted
as finding eigenfunctions of the
Hamiltonians for  the Calogero-Sutherland model for many body system, both with rational (flat case)
and hyperbolic trigonometric (bounded
case) first order differential
operators in the Hamiltonian; see e.g. \cite{Dunkl-cmp-98}
for the flat case.

The paper is organized as follows. In Section 2
we give an abbreviated introduction of weighted Bergman spaces
on bounded symmetric domains and fix some notation. In Section 3 we present the Segal-Bargmann
transform and Berezin transform on real bounded symmetric
domains, thus establish the abstract
orthogonality relation for the spherical
transforms of the Segal-Bargmann transforms
of the $L$-invariant polynomials. We identify
those
Schmid components $\mathcal P_{\n}$ which contain
non-trivial $L$-invariant polynomials in Section 4.
In Section 5 we express the radial part of the Cayley-Capelli type
operator in terms of the Dunkl operator. Their
eigenspace decomposition is done in Sections 6 and 7. In
Sections 8 and 9 we study their Segal-Bargmann
transforms and prove their orthogonality relation
and find the expansion of the Bessel and spherical
functions. We evaluate the constant
in the Plancherel formula for the symmetric space
$H/L$ in our settings in Appendix 1, and we list
all the real forms $H/L$ of a general Hermitian
symmetric space in Appendix 2.

\subsection*{Acknowledgement} 
I would like to  thank 
 Eric Opdam
for providing me a copy of the Macdonald lecture notes
\cite{Mac-leidennotes} and Charles Dunkl
for sending a copy of his works. I thank
Toshiyuki Kobayashi, Hjalmar Rosengren, Siddhartha Sahi and Harald
Upmeier for some enlighting discussions. 
The hospitality of the Newton Institute, Cambridge,
is also greatly acknowledged.

For the reader's convenience we list
the main notation used in this paper:
\begin{itemize}

\item $\mathbb D=G/K$ a complex bounded symmetric domain 
of tube type and rank $r^\prime$ in  a vector space $V_{\bc}$, $\fg=\fk+\fp$
the Cartan decomposition.

\item $V$ a real form in $V_{\bc}$ and $\tau$  the conjugation of $\Vc=V+iV$ with respect to
the real form $V$.

\item $D=H/L$ a real bounded symmetric domain in  $V$,
 $\fh=\fl+\fq$ the Cartan decomposition of $\fh$, $\fq$ is identified with $V$.

\item $r$: rank of the real bounded symmetric domain $D$, so that
the rank $r^\prime$ of the complex domain $\mathbb D$
is $r$ if $D$ is of type $D_r$,  and  it is $2r$
if $D$ is of type $C_r$.

\item $e_1, \dots, e_r$: a frame of the Jordan triple
$V$.

\item $\fa\subset \fq$ the maximal abelian
subspace of $\fq$ spanned by the vectors
$\xi_j=\xi_{e_j}$, $\Sig$ the root system
of $\fh$ with respect to $\fa$;

\item $\beta_1, \dots, \beta_r\in \fa^{\ast}$ 
the dual basis of $\frac 12 \xi_1, \dots, \frac 12\xi_r$, identified
also as linear functional on $V$.

\item $a$ the root multiplicity of $\frac{\beta_j-\beta_k}2$ and
are independent of $j, k$ and choice of the frame. Observe that
$\dim A_{jj}=1$ and $\dim B_{jj}=\iota-1$;  $\dim B_{jk}=\dim A_{jk}=a$
for type $C$ and $D_r$ ($r\ge 3$).

\item $\ga_1, \dots, \ga_r$ (if type $D$)
or $\ga_1, \ga_1^\prime \dots, \ga_{r}, \ga_r^\prime$ (if type $C$)
the Harish-Chandra strongly orthogonal
roots; $a^\prime$ the non-compact root multiplicity
of $\frac{\ga_j-\ga_k}2$.

\item $\fk= \fk^{+}_\tau + \fk^{-}_\tau
$ induced  Cartan decomposition of $\fk$
by $\tau$ with $ \fk^{+}_\tau =\fl=\fk\cap q$

\item $\ft=\ft_\tau^{-} +\ft_\tau^{+}$ the  induced decomposition of
a Cartan subalgebra $\ft$ of $\fk$.

\item $c_{\nu}$ the  normalization constant for the Berezin transform,
      $c_{\nu}^\prime$ the normalization constant
for the weighted Bergman measure,      $C_0$  the one for the Plancherel formula on $D=H/L$, $C_1$ the one
      for the integral of $L$-invariant functions on $V$ in 
      terms of polar coordinates on $\fa$.

\end{itemize}

\section{Complex bounded symmetric domains}

We recall very briefly in this and next sections
some preliminary results on bounded
symmetric domains and fix  notation;
see  \cite{Loos-bsd}
and  \cite{FK} and reference therein. 

Let $\mathbb D=G/K$ be an irreducible bounded symmetric domain
of tube type in a $d$-dimensional
complex vector space $V_{\mathbb C}=\bc^d$
of rank $r^\prime$. (The symbol $r$ will
be reserved for the rank of the real bounded symmetric
domain $D$ in Section 3.) Let $\fg=\fk+\fp$
be the Cartan decomposition
and  $\fg^{\bc}=\fp^{+}+\fk^{\bc}+\fp^-$
be the Harish-Chandra decomposition of its
complexification.
The space $V_{\mathbb C}=\fp^+$
has then a Jordan triple structure so that
the subspace $\fp$ is of the  form
\begin{equation}
\label{p-form}
\fp=\{v-Q(z)\bar v, v\in V_{\bc}\},
\end{equation}
when the elements are realized as holomorphic vector
fields, where $Q(z): \bar {V_{\bc}} \mapsto
 {V_{\bc}}$ is quadratic in $z$.
 We denote $\{x\bar y z\}$ the Jordan triple product
$ \{x\bar y z\}=(Q(x+z)- 
Q(x)- Q(z))\bar y$.
We fix
a $K$-invariant Hermitian inner product
 $(\cdot, \cdot)$ on $V_{\mathbb C}$ so that a minimal
tripotent has norm $1$. We let $dm(z)$ be
the corresponding Lebesgue measure. The Bergman
reproducing kernel up to a positive constant
is of the form $h(z, \bar w)^{-p}$ where $p$ is the genus
of $\mathbb D$, defined by $p=\frac{2d}{r^\prime}$ (for tube
domains) and $h(z, \bar w)$ is an irreducible
polynomial holomorphic in $z$ and anti-holomorphic
in $w$. 

Let $\nu>p-1$ and  consider the probability measure
$d\mu_{\nu}(z)=c_{\nu}^\prime h(z,\bar z)^{\nu -p}dm(z)$
where $c_{\nu}^\prime$
is the normalization constant,
and the corresponding weighted Bergman space $\Hnu=\Hnu(\mathbb D)$
of holomorphic functions $f$
so that
$$
\Vert f\Vert_{\nu}^2=\int_{\mathbb D}|f(z)|^2 d\mu_{\nu}(z)<\infty.
$$
It has reproducing
kernel $h(z, \bar w)^{-\nu}$.
The group $G$ acts unitarily on
$\Hnu$ via the following
\begin{equation}
\label{pi-nu}
\pi_{\nu} f(z)= J_{g^{-1}}(z)^{\frac \nu p}f(g^{-1}z),  
\end{equation}
and it forms a unitary projective representation of $G$.

Let $\mathcal F_{\nu}=\mathcal F_{\nu} (V_{\bc})$, for $\nu >0$,
be the Fock space of entire function on $V_{\bc}$
with the norm
defined by
$$
(p, q)_{\mathcal F_\nu}=\frac{\nu^d}{\pi^d}
\int_{V_{\bc}}p(z)\overline{q(z)}e^{-\nu (z, z)} dm(z).
$$
Thus it has reproducing kernel $e^{\nu(z, w)}$.
When $\nu=1$ we write for simplicity $\mathcal F=\mathcal F_{1}$.
The norm in $\mathcal F$ can alternatively
defined by
\begin{equation}
  \label{fock-nu}
(f, g)_{\mathcal F_\nu}=\nu^{-\deg (p)}
f(\partial_z) g^\ast(z)\big{\vert}_{z=0}  
\end{equation}
for polynomials $f$ and $g$, where $g^\ast$ is obtained
from $g$ by taking formerly the complex conjugate of
the coefficients of the monomials
(in terms of an orthonormal basis in $V_{\bc}$).

The first part  of the next result
is due to Hua \cite{Hua} for classical domains and
Schmid \cite{Schmid} for general domains, the second and third part
is due to Faraut and Koranyi \cite{FK}, Theorem 3.8; this
result will be of fundamental importance for our work here.
Fix a Cartan subalgebra of $\fk$ and
let $\ga_1>\dots>\ga_{r^\prime}$ be 
the Harish-Chandra strongly orthogonal
roots. Let $a^\prime$ be the
root multiplicity of $\frac{\ga_j-\ga_k}2$
in $\fp^+$.

\begin{theo+} The space $  \mathcal P$
of holomorphic polynomials on $V_{\bc}$ decomposes into irreducible
subspaces under $Ad(K)$, with multiplicity one as: 
\begin{equation}
\label{hua-sch}
  \mathcal P\cong \sum_{\underline{\mathbf{n}} \ge 0} \Pn.
\end{equation}
Each $\Pn$ is of lowest weight $-\underline{\mathbf{n}}%
=-(n_1\gamma_1 +\cdots +n_r \gamma_r)$ with $n_1\ge \dots \ge n_r\ge 0$.
The quotient of the norms of
 a polynomial $f\in \Pn$ 
in the Fock space $\mathcal F$ and in the weighted
Bergman space $\Hnu$ is given by
$$
\frac{\Vert f\Vert_{\mathcal F}}{
\Vert f\Vert_{\Hnu}
}=\sqrt{(\nu)_{\n}}
$$
where
$$
(\nu)_{\n}=\prod_{j=1}^{r^\prime} (v-\frac {a^\prime}2(j-1))_{n_j}=
\prod_{j=1}^r \prod_{k=1}^{n_j}
(v-\frac{a^\prime}2(j-1) +k-1).
$$
is the generalized Pochhammer symbol.
\end{theo+}

Let $K_{\n}$ be the reproducing kernel
of $\Pn$ with in the the Fock space
$\Fnu$ for $\nu=1$. Then as a consequence we have,
writing $|\n|=n_1+\dots +n_{r^\prime}$,
\begin{equation}
  \label{exp-gau}
  e^{\nu(z, w)} =\sum_{\n} v^{|\n|}K_{\n}(z, w) 
\end{equation}
and 
\begin{equation}
  \label{exp-h}
  h(z, \bar w)^{-\nu} =\sum_{\n}(\nu)_{\n}K_{\n}(z, \bar w).
\end{equation}

\section{Berezin and Bargmann transform
on real bounded symmetric domains}

 Let $V$ be a real form of $\Vc$ and $\tau$
 the complex conjugation with respect to $V$. Suppose $\tau(\mathbb D)=\mathbb D$,
namely,
$\tau$ fixes the bounded symmetric domain. 
Then the real form
$D=\mathbb D\cap V$ is called a \textit{real bounded symmetric domain}.
In this case the triple product $D(x, \bar y) z=
\{x\bar y z\}$ restricted on $V$
defines also a triple product on $V$.

A complete
list of real bounded symmetric domains $D$
is given in \cite{Loos-bsd}. 
As a Riemannian symmetric space, $D=H/L$, where $H$ is
the connected component of
the subgroup of $G$ of biholomorphic transformations
 of $\mathbb D$ which keep $D$  invariant. 
The coset space $G/H$  is
 called a causal symmetric space,
a complete list of the pairs $(G, H)$ can be found
in  e.g. \cite{Ola-dg} and \cite{Ola-Hilg}, see also Appendix 2.

The $H$-invariant measure on $D$ is 
\begin{equation}
  \label{H-invmea}
{h(z, \bar z)^{-\frac p 2}}dm(z),  
\end{equation}
and  $H$ acts unitarily on $L^2(D, h^{-\frac p 2}dm)$ via change of 
variables,
\begin{equation}\label{pi-0}
  \pi_0(g)f(x)=f(g^{-1}x), \quad g\in H.
\end{equation}

We describe
briefly some  algebraic and geometric structures of the domain $D$.

 Let $\fh=\fq\oplus\fl $
be the Cartan decomposition of the Lie algebra $\fh$ of $H$.
Similar to (\ref{p-form})
we have
\begin{equation}
\label{q-form}
\fq=\{v-Q(z)\bar v, v\in V\}.
\end{equation}
We thus identify $\fq$ with the underlying space $V$
via the mapping $v\mapsto \xi_v(z)$.

Let  $r$  the rank of $D=H/L$, so that
 the rank $r^\prime$ of $\mathbb D$ will be $r$ or $2r$; see below.
Let $\{e_j, j=1, \dots, r\}$ be a frame   of minimal tripotents in $V$.
The  elements $\xi_j=\xi_{e_j}=e_j-Q(z)\bar {e_j}$, $j=1, \dots, r$,
span  a maximal subspace $\fa$ 
of $\fq$ of dimension $r$. Let $\{\b_j, j=1, \dots, r\}$ in
$\fa^\ast$ be
the dual basis of $\frac 12 \xi_j$,
$$
\beta_j(\xi_k)=2\delta_{j, k},
$$
where $\delta_{j, k}$ 
is the Kronecker symbol. Then the root
system  $\Sig(\fg, \fa)$ 
is of types $A_r$, $C_r$ or $D_r$. Type $A_r$
corresponds to the case that $V$ is a formal real Jordan
algebra; we will be only concerned with type $C_r$
and type $D_r (r\ge 3)$:
 $$\Sig(\fg, \fa)=\{\pm \beta_j, \frac{\beta_j\pm \beta_k}2\} ;$$
and respectively
 $$
\Sig(\fg, \fa)=\{\frac{\beta_j\pm \b_k}2\}. $$ 
We let,  as in \cite{gkz-bere-rbsd},
$a$ be the root multiplicity of $\frac {\b_j\pm \b_k}2$
and $\iota-1$ that of $\b_j$.
Arrange an ordering  of the roots so that
$$
\b_1<\b_2<\cdots <\b_r.
$$
The ranks  $r^\prime$ and $r$ and the multiplicities  
$a^\prime$, $a$ and $\iota-1$ are related to
each other:
\begin{equation}
  \label{r-r-D}
r^\prime=r, a^\prime=2a  
\end{equation}
if it is type $D_r(r\ge 3)$, and
\begin{equation}
  \label{r-r-C}
r^\prime=2r, a=2a^\prime, \iota=2+a^\prime;
\end{equation}
if it is type $C_r$ for $r\ge 2$; $a^\prime =0$
and $\iota=2+a^\prime$ if $r=1$.
Those relation
can easily be obtained by calculations
of the dimensions of the root spaces
or by a case by case check of the
table in the Appendix 2.

We note that with the above normalization
of the inner product on $\Vc$
a  minimal tripotent of $V$ has norm ${\sqrt 2}$  if
the root system is of Type C, otherwise
it is of norm $1$; see below. (The inner product in
\cite{gkz-bere-rbsd} is normalized so that a minimal
tripotent of $V$ always has norm $1$.)

To consider the branching law of
$(\mathcal H_{\nu}, \pi_{\nu})$ of
$G$ under $H$ we let $R$ be the restriction map 
(\cite{oo-restric-1995})
$R: H_{\nu} \to C^\infty (D)$ by
\begin{equation}\label{def-R}
Rf(x)=f(x)h(x, \bar x)^{\frac{\nu}2}, \quad x\in D.
\end{equation}
Then $R$ is an $H$-intertwining map, as one can easily checks from
the transformation properties of $h(x, \bar x)$. Consider its formal
conjugate operator $R^\ast$
from $L^2(D, d\mu_0)$ to $H^\nu$ and form 
the operator $R^\ast R$
on $L^2(D, \ome)$. It is 
\begin{equation}
  \label{RRstar}
RR^\ast=\frac 1{c_\nu}B_\nu,
\end{equation}
where
\begin{equation}
  \label{def-Bnu}
B_\nu f(z)=c_\nu\int_Df(w)\frac{h(z, \bar )^\frac \nu 2 h(w, \bar w)^\frac \nu 2}{
h(z, \bar w)^{\nu}} \frac {dm(w)}{h(w, \bar w)^{\frac p 2}}.  
\end{equation}
is the (normalized) Berezin transform
and the constant $c_\nu$ is  such that
$B_\nu 1 =1$.  The constant $c_\nu$ is
evaluated in \cite{gkz-bere-rbsd}.
Let 
\begin{equation}
\label{polar-deco}
  R=|R|U
\end{equation}
be the polar decomposition of $R$. Thus
$|R|^2=RR^\ast$.

It is proved in \cite{gkz-bere-rbsd}
that when $\nu >p-1$ for
type C and suppose $\nu >\frac{p}2-1$
the  Berezin transform $RR^\ast$
is bounded, and that the multiplier
$h^{\frac{\nu}2}$ is in $L^1(D, h^{-\frac p2}dm(z))$.
Since for those values of $\nu$
the space $\Hnu$
contains all polynomials, thus
the range of $R$ contains the functions
of the form $h(z, \bar z)^{\frac \nu 2}p(z, \bar z)$
where $p(z, \bar z)$ are polynomials
of $z$ and it clearly forms a dense subspace
in $L^2(D, h^{-\frac p2}dm(z))$.
This proves 

\begin{prop+}Suppose $\nu >p-1$ for
type C and suppose $\nu >\frac{p}2-1$ for
type D.
 The operator $U$ is  unitary
and intertwines $H$-actions (\ref{pi-nu})
on $\Hnu$ and (\ref{pi-0})
 onto $L^2(D)$. 
\end{prop+}

\begin{defi+} The unitary operator $U$ is called
a (generalized) Segal-Bargmann transform.
\end{defi+}

Let $\phi_{\blam}$ be the spherical function on 
$D=H/L$. Let
$$
\widehat f(\blam)=
\int_D f(z)\phi_{\blam}(z)\frac{dm(z)}{h(z, \bar z)
^{\frac p2}}
$$
be the spherical transform
on $D$, where $f$ is a $L$-invariant
$C^\infty$-function on $D$ with compact support. 
It is well known (\cite{Hel1}, Chapter IV)
that the map extends to
a unitary operator from $L^2(D)^L$ onto 
$L^2(\fa^\ast, C_o|c(\blam)|^{-2}d\blam)^W$ of 
$W$-invariant functions in the $L^2$ space on
$\fa^\ast$. Here $c(\blam)$ is the
Harish-Chandra $c$-function, with
the same normalization
as in  \cite{Hel1}, Chapter IV,
$d\blam$ is the regularly normalized
measure and the constant $C_0$ can be
evaluated by using the evaluation of the constant
$c_{\nu}$ in \cite{gkz-bere-rbsd} and
the evaluation of the Selberg integral, see Appendix 1.

\begin{coro+}\label{realize!} 
The map $f\mapsto \widehat{Uf}$
is a unitary operator from $\Hnu^L$
onto $L^2(\fa^\ast,  C_o|c(\blam)|^{-2}d\blam)^W$.
\end{coro+}

In particular if we find a canonical orthogonal
basis $\Hnu^L$, their image under
the above unitary operator then gives
an orthogonal basis $L^2(\fa^\ast,  C_o|c(\blam)|^{-2}d\blam)^W$.

\section{Existence  of the $L$-invariant holomorphic polynomials}

We identify now those polynomial spaces $\Pn$ that
contain $L$-invariant vectors.

\begin{lemm+} In the decomposition (\ref{hua-sch})
the component $(\mathcal P_{\n})^L\ne 0$ if and only if
\begin{equation}
\label{n-m-rel-c}
\n=(m_1, m_1, m_2, m_2,\dots, m_r, m_r)=\sum_{j=1}^rm_j(\ga_{2j-1} +\ga_{2j})
\end{equation}
if $\Sig$ is type $C_r$,
and
\begin{equation}
\label{n-m-rel-d}
\n=(2m_1, 2m_2, \dots, 2 m_r)+m(1, 1, \dots 1)=
\sum_{j=1}^r(2m_j+m)\ga_{j} 
\end{equation}
if $\Sig$ is of type $D_r$,
where in all cases $m_j$ and $m$ are nonnegative integers
and $m_1\ge m_2\ge \dots \ge m_r\ge 0$.
\end{lemm+}

\begin{proof}  The involution $\tau$
on $\Vc$ induced an involution on $\fk$
and its fixed point set is $\fl$, thus
$(\fk, \fl)$ is a symmetric pair.
Let $\fk=\fl+\fk_{\tau}^{-}$
be the corresponding Cartan decomposition. Let
 $\ft_{\tau}^-\subset \fk_{\tau}^{-}$
be a maximal abelian subspace of $\fk_{\tau}^{-}$,
and $\ft=\fk_{\tau}^{+}
+\fk_{\tau}^{-}$ a Cartan subalgebra of $\fk$
with $\ft_{\tau}^+\subset \fl$.
 We
apply  the Cartan-Helgason theorem
\cite{Hel2}, Chapter V, Theorem 4.1, 
to identify the $L$-spherical representations,
which asserts in our case that the representation $\Pn$
contains a $L$-fixed vector if and only if
\begin{equation}
  \label{van-con}
\n\big{\vert}_{(\fk_{\tau}^{+})^{\bc}
} =0 
\end{equation}
and 
\begin{equation}
  \label{int-con}
\frac{(\n, \a)}{(\a, \a)}\in \mathbb N,
\end{equation}
for all roots
$\a$ in the root system $\Sig(\fk^{\bc}, (\ft_{\tau}^{-})^{\mathbb C})$
of $\fk^{\bc}$ with respect
to $\ft_{\tau}^{-}$, in which case
$\Pn^L$ is one-dimensional. Here $\mathbb N$ is the set of
nonnegative integers.

Consider the type $D_r$ first. By (3.1) the frame
$\{e_1, e_2, \dots, e_r\}$
of minimal tripotents in $V$ is also
a frame in $V_{\bc}$. Recall the involution $\tau$
induced by the complex conjugation on   $V_{\bc}$
with respect to  $V$. We have $\tau(\xi_{e_j})=\xi_{e_j}$
$\tau(\xi_{ie_j})=-\xi_{e_j}$. Thus
$$
\tau(iD(e_j, e_j))=-iD(e_j, e_j),
$$
since
$$
iD(e_j, e_j)=\frac 12[\xi_{e_j}, \xi_{ie_j} ],
$$
so that  $iD(e_j, e_j)\in \ft_{\tau}^-$, $j=1, \dots, r$,
and they span an $r$-dimensional abelian subspace
$i(\br D(e_1, e_1) +\dots +\br D(e_1, e_1))$
of $\ft_{\tau}^-$, whereas the dimension of $\ft_{\tau}^-$
is $\dim(\ft_{\tau}^-)=\text{rank}(K/L)$.
 We claim that $\text{rank}(K/L)>r$.
Indeed the  symmetric pair $(\fk, \fl)$
is $( s(\mathfrak{u}(r)\oplus \mathfrak{u}(r)),
\mathfrak{so}(r)\oplus \mathfrak{so}(r))$,
$(\mathfrak{u}(2r), \mathfrak{so}(2r))$
 or $ (\mathfrak{so}(10)+\mathbb R ,\mathfrak{sp}(4))$
(with $r=3$), the rank of the first pair
being $2r-1 >r$, the second $2r>r$ and
the third $7>3$ (see \cite{Hel1}, Table V, p. 518). This means that
the subspace $i(\br D(e_1, e_1) +\dots +\br D(e_1, e_1))$
is a nontrivial subspace of $\fk_{\tau}^-$.
So the vanishing
condition (\ref{van-con}) is already satisfied
since $\n$ is vanishing on
the orthogonal complement
of  $i(\br D(e_1, e_1) +\dots +\br D(e_1, e_1))$
in $\ft$. To check the
second condition, recall (see \cite{Schmid}, formulas (16)
and (17))
that all roots of in $\Sig(\fk^{\bc}, (\ft_{\tau}^{-})^{\mathbb C})$
are of the form $\frac{\ga_j-\ga_k}2 + \a$, with $\a$ orthogonal to $\frac{\ga_j-\ga_k}2$ and if nonzero,
 $\Vert\a\Vert^2=\Vert\frac{\ga_j-\ga_k}2\Vert^2$.
Notice that since $\frac{\ga_j-\ga_k}2$,
viewed as linear functional on $\ft_\tau^-$, only
span a $r-1$-dimensional subspace, thus there exist roots
in $\Sig(\fk^{\bc}, (\ft_{\tau}^{-})^{\mathbb C})$
of the form $\frac{\ga_j-\ga_k}2 + \a$
with nonzero $\a$ for all $j,k=1, \dots, r$, $j\ne k$.
Thus by the Cartan-Helgason theorem, 
the polynomial representation $\n$
is $L$-spherical if and only if
$$
\frac{\langle \n, \frac{\ga_j-\ga_k}2 + \a\rangle}
{\langle \frac{\ga_j-\ga_k}2 + \a, \frac
{\ga_j-\ga_k}2 + \a\rangle}=
\frac 12(n_j-n_k)\in \mathbb Z_{\ge 0},
$$
for all $j<k$, 
and
$$
\frac{\langle \n, \frac{\ga_j-\ga_k}2\rangle}
{\langle \frac{\ga_j-\ga_k}2, \frac
{\ga_j-\ga_k}2 \rangle}=
(n_j-n_k)\in \mathbb Z_{\ge 0},
$$
when $\a=0$. The first condition clearly
implies the second and it is just our stated
condition.

Now consider type $C_r$. 
There
are only three cases and we study them case by case.
Consider the case
 $(\fk, \fl)=(s(\mathfrak u(2r) \oplus
\mathfrak u(2r)), \mathfrak {sp}(r) \oplus
 \mathfrak {sp}(r)
  )$. The highest weight
$\n=(m_1\varepsi_1+m_2\varepsi_2
 +\dots +m_{2r}\varepsi_{2r})
\otimes
(m_1\varepsi_1+m_2\varepsi_2
 +\dots +m_{2r}\varepsi_{2r})^\ast
$
where $\varepsi_j$ is the dual of diagonal matrix
with $j$th entry being $1$ and rest $0$,
in the standard matrix
representation of $\mathfrak{u}(2r)$,
and   $(m_1\varepsi_1+m_2\varepsi_2
 +\dots +m_{2r}\varepsi_{2r})^\ast$
is the contra-gradient representation.
The representation $\m$ has
a $\fl=\mathfrak {sp}(r) \oplus
 \mathfrak {sp}(r)$-fixed vector if and only if
the representation $m_1\varepsi_1+m_2\varepsi_2
 +\dots +m_{2r}\varepsi_{2r}$ has
a $\mathfrak {sp}(r)$-fixed vector, and the later
happens precisely when $m_1=m_2, m_3=m_4, \dots, m_{2r-1}=m_{2r}$, namely our condition,
 again by the Cartan-Helgason theorem, since
$ (\mathfrak {su}(2r),  \mathfrak {sp}(r))$
is a symmetric pair. Now let
 $(\fk, \fl)=(\mathfrak u(2r), \mathfrak {sp}(r))$. The Harish-Chandra roots are $\ga_j=2\varepsi_{2j-1}, \ga_j^\prime =2\varepsi_{2j}$. Our result immediately follows 
from 
the above argument.
Finally
consider 
 $(\fk, \fl)=(\mathfrak {so}(2)\oplus
 \mathfrak {so}(p), \mathfrak {so}(p))$. The highest
weight $\n=m_2(\ga_1+\ga_2) +(m_1-m_2)(\ga_1)$,
and as representations $\n=m_2(\ga_1+\ga_2)\otimes
(m_1-m_2)\ga_1 $ with $m_2(\ga_1+\ga_2)$ being the
trivial representation of $\mathfrak {so}(p)$
and $(m_1-m_2)\ga_1 $
the representation of $\mathfrak {so}(p)$
on the spherical harmonics of degree $m_1-m_2$, the later
has an $\mathfrak {so}(p)$-invariant vector if
and only if $m_1-m_2$ is even.
This completes the proof.
\end{proof}

\begin{rema+}
In \cite{Kraemer} Kr\"amer has
given a classification of all spherical
connected subgroups $L$ of
 a simple compact Lie group $K$ (Tabelle 1, loc. cit)
and listed all the fundamental spherical
highest weights.
Our result can also be deduced from that
list. One may also prove the above
result somewhat by slightly general argument
by following the classification
as in \cite{Ola-dg} and \cite{Loos-bsd}.
 Moreover the above result
for type C holds also for  type BC
and that for type D holds also for type B (with $m=0$).
\end{rema+}

Let $\Del(z)$
be the determinant polynomial of
the Jordan triple  $V_{\bc}$ and $\Del(\partial)$
the corresponding differential operator which
we call the Cayley type operator.
The next result gives
the eigevalue of the operator
$\Del(z)^{-\a}\Del(\partial)\Del(z)^{\a +1}$,
sometimes also named as the Cayley-Capelli  operator,
under the Schmid decomposition (\ref{hua-sch}),
and its follows  easily by using Theorem 2.1;
 see e.g. \cite{Yan-inv-eig}, 
\cite{Wallach-sha}.

\begin{lemm+}\label{cayley-dia}
 Let $r^\prime$ be the rank of $\mathbb D$.
The differential operators $\Del(z)^{-\a}\Del(\partial)\Del(z)^{\a +1}$
for $r^\prime$ different nonnegative integers $\a$ form
a system of generators of $K$-invariant
differential operators on $ \mathcal P(V_{\mathbb C})$.
A holomorphic polynomial $f\in \mathcal P(V_{\mathbb C})$
is in the space $\Pn$ if and only if it is a solution of
the system of differential equations
\begin{equation}
\label{cayley}
\Del(z)^{-\a}\Del(\partial)\Del(z)^{\a +1}p(z)
=\prod_{j=1}^{r^\prime}(\frac {a^\prime}2(r^\prime-j)+1+\a +n_j)p(z),
\end{equation}
for $r^\prime$ different nonnegative integers $\a$.
\end{lemm+} 

In the next sections we will find
the $L$-invariant polynomials in $\Pn$.

\section{Determination  of the $L$-invariant holomorphic polynomials: Some general results }

We consider the Chevalley restriction
map $\Res$ from $\PV^L$ onto the space $\Pa^W$ of
 $W$-invariant polynomials on $\fa$.  By using
the Dunkl operator
we define an inner product on $\Pa^W$
and we prove that the restriction map  $\Res$  is
an unitary map from $\PV^L$ with the Fock space norm
onto  $\Pa^W$. This reduces the problem of finding 
and calculating of $L$-invariant
polynomials to the corresponding one of 
$W$-invariant polynomials on $\fa$.

Some consideration that follows will be true
for some general root system, we shall,
however, only consider the 
root systems $\Sig$ as in Section 3.

We fix  a frame  $\{e_j, j=1, \dots, r\}$ 
 of minimal tripotent in $V$,
enumerated so that $e_j\in V\subset V^{\bc}=\fp^+$
is a root vector of the Harish-Chandra
orthogonal root $\ga_j$, $j=1, \dots r$, if $\Sig$
is of type C, and $e_j$ is a sum of
two root vectors with roots $\ga_{2j-1}$ and $\ga_{2j}$
for type D.

For any root $\a\in \Sig$
let $r_{\a}\in W$ be the  reflection
defined by $\a$.
We recall the Dunkl operator
\cite{Dunkl-tams},
$$
D_j=\partial_j +\frac 12 \sum_{\a\in \Sig^+}m_{\a}\frac{\a(\xi_j)}{\a(x)}
(1-r_{\a})
$$
acting on polynomials $f(x)$ on $\fa$, where
the $r_{\a}$ acts $f(x)$ via
$$
(r_{\a}f)(x)=f(r_{\a}^{-1}x).
$$
The operators $D_j, j=1, \dots, r$ are pairwise commuting,
and thus define an isomorphism between the ring of polynomials
on $\fa$ and the ring of difference-differential operators
generated by them. We can thus define, for any
polynomial $x=x_1e_1+\dots +x_r e_r\to f(x)=f(x_1, \dots, x_r)$ 
the operator $f(D)$ by assigning $D_j$ to
the polynomials $x_j$, $j=1, \dots, r$.

\begin{defi+}(Dunkl \cite{Dunkl-cjm}) Let $\Sig$
be the root system of type $C$ or type $D$.
The $\Sig$-inner product on $(\Pa)^W$ is defined
by
$$
(f, g)_{\Sig}=f(D_x)g^{\ast}(x)\big{\vert}_{x=0},
$$
if $\Sig$ is of type D, and
$$
(f, g)_{\Sig}=f(\frac 12 D_x)g^{\ast} (x)\big{\vert}_{x=0}
$$
if $\Sig$ is of type C.
\end{defi+}

The discrepancy for type C here
is due to the unmatched norms of minimal tripotents
in $V_{\bc}$ and $V$.

\begin{prop+}The restriction map $\Res$
is an isometric mapping from $(\PV^L, \Vert\cdot\Vert_{\mathcal F})$
onto  $(\Pa^W, \Vert\cdot\Vert_{\Sig})$.
\end{prop+}
\begin{proof} We consider the operators
  \begin{equation}
    \label{E}
E=\frac 12 (z_1^2+\dots +z_d^2),    
  \end{equation}
  \begin{equation}
    \label{F}
F=-\frac 12 (\partial_1^2+\dots +\partial_d^2),
  \end{equation}
and
  \begin{equation}
    \label{H}
    H=(z_1\partial_1+\dots +z_d\partial_d) +\frac d2
  \end{equation}
acting on the space $\PV^L$. Then it is clear that
they form the Lie algebra $\mathfrak{sl}(2, \bc)$:
$$
[E, F]=H, \quad [H, E]=2E, \quad [H, F]=2F.
$$
Let us, following
Heckman \cite{Heckman-remark}, define similarly
  \begin{equation}
E_0=\frac 12 (x_1^2+\dots +x_r^2),    
  \end{equation}
  \begin{equation}
F_0=-\frac 12 (D_1^2+\dots +D_r^2),
  \end{equation}
for type D, and
  \begin{equation}
E_0=x_1^2+\dots +x_r^2,    
  \end{equation}
  \begin{equation}
F_0=-\frac 14 (D_1^2+\dots +D_r^2),
  \end{equation}
for type C (again due to
the ill-matching of the norms of minimal tripotents
in $V_{\bc}$ and $V$),
 and
  \begin{equation}
    H_0= (x_1\partial_1+\dots +x_r\partial_d) +\frac 12(r + \sum_{\a\in R}m_\a),
  \end{equation}
for all types. They form a copy of the Lie algebra $\mathfrak{sl}(2, 
\mathbb C)$.
We claim that
$$
\Res E=E_0, \quad  \Res F=F_0, \quad \Res H=H_0.
$$
The second is trivial, the third follows since
 that the dimension 
$d=\dim_{\bc }V_{\bc}=\dim_{\br}\fq=\frac r2 +\frac
12 \sum_{\a\in R}m_\a$. The first is just
the formula for the radial part of the Laplace operator
on $\fq$; see \cite{Hel2}, Proposition 3.13 for
the formula and \cite{Dunkl-tams} (or \cite{Heckman-remark})
for the
calculation of $F_0$. It is proved in \cite{Heckman-remark} Proposition
3.4 that for any polynomial $p$ of degree $m$ on $\fa$,
viewed as multiplication operator on $\Pa$, namely
in $Aut(\Pa)$,
\begin{equation}
  \label{pD}
p(D)=(-1)^m\frac 1{m!}\ad(F_0)^m (p),
\end{equation}
where the Lie algebra $\mathfrak{sl}(2, \bc)$ is acting on
$Aut(\Pa)$ via the adjoint action.
The essentially same (even easier) calculation
shows that, for polynomial $P$ on $\PV$ we have
\begin{equation}
  \label{PD}
P(\partial)=(-1)^m\frac 1{m!}\ad(F)^m (P)  
\end{equation}
For any $P, Q\in \PV^L$, 
let $p=\Res P$, $q=\Res Q$; if $x\in \fa$,
\begin{equation*}
\begin{split}
P(\partial)Q^\ast(x)
&=(\Res  P(\partial)Q^\ast)(x)\\
&= \left(\Res ((-1)^m\frac 1{m!}\ad(F)^m (P)
Q^\ast)\right) (x)\\
& =
\left((-1)^m\frac 1{m!}\ad(F_0)^m (\Res P) \Res Q^\ast\right)(x)\\
&=p(\partial)q^\ast(x)
\end{split}
\end{equation*}
and
$$
(P, Q)_{\mathcal F}=
P(\partial)Q^\ast(0)
=p(\partial)q^\ast(0)
=(p, q)_{\Sig}
=(\Res P, \Res Q)_{\Sig},
$$
completing the proof.
\end{proof}

It is noted  in \cite{Heckman-remark} that the above idea
of reducing computation to second order
operators goes back to Harish-Chandra.

\begin{rema+}The above result clarifies the significance
of the Dunkl operator and the inner product $(\cdot, \cdot)_{\Sig}$.
Moreover it gives an isometric version of the Chevalley
restriction theorem (\cite{Hel2}, Chapter 2, Corollary 5.12). It
seems to the author that this was not been known before. 
\end{rema+}

As a consequence we get
\begin{coro+} The Cayley  operator $\Res \Del(\partial)$
on $L$-invariant polynomials is given by
$$
(\Res \Del(\partial))f(x)=\prod_{j=1}^r D_j f(x)
$$
for type $D$
and 
$$
(\Res \Del(\partial))f(x)=2^{-2r}\prod_{j=1}^r D_j^2 f(x)
$$
for type $C$.
\end{coro+}
\begin{proof} We consider only the
 type $D_r$, the type $C_r$-case can be
proved similarly. Recall that the determinant
function
$\Del$  on $V_{\bc}$ is of
degree $r$ and  when restricted on $\fa$ it is
$$
\Res \Del(x)=\Del(x_1e_1 +\dots +x_re_r)=
\prod_{j=1}^r x_j=(\prod_{j=1}^r \frac{\b_j}2)(x).
$$
We calculate now the adjoint operator
of $\Del$ and $\Res\Del$. The adjoint operator
of multiplication by $\Del(z)$ on $\PV$ with
respect to the Fock norm is
 $\Del(\partial)$, $\Del^\ast=\Del(\partial)$ and respectively
the multiplication by $\prod_{j=1}^r \frac{\b_j}2$
is $(\prod_{j=1}^r \frac {\b_j}2)^\ast=\prod_{j=1}^r D_j$
with respect to the $\Sig$-norm by the definition of the inner product. Thus
for $P, Q\in \PV^L$ we have, by the preceding
proposition, 
\begin{equation*}
\begin{split}
&\quad (\Res (\Del(\partial) P), \Res Q)_R=
(\Del(\partial) P,  Q)_{\mathcal F}=
(P, \Del Q)_{\mathcal F}\\
&=(\Res P, \Res \Del Q)_{\Sig}=
(\Res P, (\prod_{j=1}^r\frac{\b_j}2)\Res \Del Q)_{\Sig}\\
&=
((\prod_{j=1}^r\frac {\b_j}2)^\ast\Res P, \Res Q)_{\Sig}=
((\prod_{j=1}^rD_j)\Res P, \Res Q)_{\Sig}
\end{split}
\end{equation*}
proving the result.
\end{proof}

\section{$L$-invariant holomorphic polynomials: Type $C_r$}

In this section we find the $L$-invariant
polynomials in $\Pn^L$ and calculate
their norm in the Fock space.  We will express 
them in terms of the Jack symmetric
polynomials. For that purpose
we recall some basic facts.

Let 
\begin{equation}
\label{duopa}
D_j=D_j^A=\partial_j + \frac a2 \sum_{i\ne j}\frac{1}{y_j-y_i)}(1-s_{ij})  
\end{equation}
be the Dunkl operator (\cite{Dunkl-tams},
\cite{Heckman-remark}) acting on functions
on the $r$-dimensional vector space $\br^r$, where
 the superscript indicates
that the underlying root system is of type A. Let
\begin{equation}
\label{chopa}
U_j=U_j^{A}=D_j y_j-\frac a2 \sum_{i<j}s_{ij}
\end{equation}
be the Cherednik operator for type $A$
(\cite{Dunkl-cmp-98} and \cite{Chered-uni}).
Then $U_j^{A}$, $j=1, \dots, r$ are commuting
operators. The relation between $U_j^{A}$ and
the Cayley type operator $\prod_{j=1}^rD_j
\prod_{j=1}^r y_j$ is the following. (It can be viewed
as a symmetric invariant analogue
of the Cayley-Capelli identity, which express
the product $\Del(\partial)\Del(z)$ as 
another determinant.)

\begin{lemm+} The following identity hold
$$
(\prod_{j=1}^r D_j)(\prod_{j=1}^r y_j)
=\prod_{j=1}^r U_j
$$
\end{lemm+}
\begin{proof} The proof relies on
the commutation relations $D_i s_{ij}=s_{ij}D_j$
and $[D_i, x_j]=-\frac a2 s_{ij}$. From this
we deduce that
\begin{equation*}
\begin{split}
\prod_{j=1}^r D_j\prod_{j=1}^r y_j
&=D_1\dots D_{r-1} D_r y_1 y_2\dots y_r\\
&=D_1\dots D_{r-1} y_1 D_r y_2\dots y_r-
\frac a2 (\prod_{j=1}^{r-1} D_j)(\prod_{j=1}^{r-1} y_j)s_{1r};
\end{split}
\end{equation*}
repeating the argument (moving $D_r$ until it reaches $y_r$) we get
\begin{equation}
\prod_{j=1}^r D_j\prod_{j=1}^r y_j
=D_1\dots D_{r-1}  y_1 y_2\dots y_{r-1} U_r.
\end{equation}
Performing the above computation with $D_{r-1}$ and so on
proves the formula.
\end{proof}

Essentially the same computation gives
\begin{equation}
  \label{Cayley-A}
(\prod_{j=1}^r y_j)^{-\a}\prod_{j=1}^rD_j (\prod_{j=1}^r y_j)^{1+\a}=
\prod_{j=1}^r(U_j+\a).
\end{equation}

\begin{rema+} The Cherednik operators $U_j^{A}$ are introduced
before (see \cite{Dunkl-cmp-98} for root systems
of type B and reference therein)
in order to study the non-symmetric Jack polynomials.
The above lemma shows that it can also be obtained,
though less systematically, in trying to write the Cayley-type operator
$(\prod_{j=1}^r D_j)(\prod_{j=1}^r x_j)$ as a product
of $r$ operators.
\end{rema+}

The Jack symmetric polynomial
are then a sum of joint (nonsymmetric) eigenfunctions of $U_j$. Let
$\Ome_{\m}$ be the normalized Jack symmetric polynomial, so
that it is a symmetric
eigenpolynomial of the operators $\prod_{j=1}^r(U_j+\a)$
\begin{equation}
  \label{jack-def}
\prod_{j=1}^r(U_j^{A}+\a)
\Ome_{\m}=
(\prod_{j=1}^{r}(\frac a2(r-k)+1+\a +m_k))\Ome_{\m} (y_1, \dots, y_r)
\end{equation}
for all nonnegative integers $\a$, normalized so that
$$
\Ome_{\m}(1, \dots, 1)=1.
$$

Let $D_j$ be the Dunkl operator
on $\fa=\br e_1 +\dots \br e_r$
for the root system $\Sig$ of type C,
$$
D_j=\partial_j +\frac{\iota-1}2
\frac{1}{x_j}(1-s_{j})+\frac a2\sum_{i\ne j}
(\frac{1}{x_j -x_i}(1-s_{ij})+\frac{1}{x_j +x_i}(1-\sig_{ij}))$$
where $\sig_j$, $s_{ij}$ and $\sig_{ij}$
are the reflections in Weyl group 
corresponding to the roots $\ga_j$, $\frac 12(\ga_i-\ga_j)$
and respectively  $\frac 12(\ga_i+\ga_j)$.

Recall that the Hermitian form on $V$ is normalized
so that $e_j$ has norm $\sqrt 2$.
\begin{lemm+} The restriction of the operator
$\Res(\Del(\partial_x)\Del(x))$ on $\fa$ is given by
\begin{equation*}
\begin{split}
&\quad \Res(\Del(\partial_x)\Del(x))\\
&=2^{-2r}\prod_{j=1}^r D_j^2 \prod_{j=1}^r x_j^2 
=2^{-2r}\prod_{j=1}^r D_j \prod_{j=1}^r
(D_jx_j -\frac a2 \sum_{i<j}(\sig_{ij} +s_{ij}))
\prod_{j=1}^r x_j\\
&=2^{-2r}\prod_{j=1}^r\left(D_jx_j -\frac a2 \sum_{i<j}(\sig_{ij} +s_{ij}) +1-
(\iota-1)\sig_j)\right)
\prod_{j=1}^r\left(D_jx_j -\frac a2 \sum_{i<j}(\sig_{ij} +s_{ij})\right)
\end{split}
\end{equation*}
\end{lemm+}
\begin{proof} The first equality follows by
 Corollary 5.4, for the restriction
of $\Del(x)$ on $\fa$ is $\prod_{j=1}^r x_j^2$.
The rest of the proof is similar to that of Lemma 6.1.
We have the product formula
$$
(\prod_{j=1}^rD_j)(\prod_{j=1}^rx_j)=
\prod_{j=1}^r(D_jx_j -\frac a2 \sum_{i<j}(\sig_{ij} +s_{ij}));
$$
and shift formula
$$
(\prod_{k=1}^rD_k)(D_jx_j -\frac a2 \sum_{i<j}(\sig_{ij} +s_{ij}))
=(D_jx_j -\frac a2 \sum_{i<j}(\sig_{ij} +s_{ij})+1-(\iota-1)\sig_j)
(\prod_{k=1}^rD_k),
$$
which can be obtained  by repeatedly using the commutation relation
$$
[D_j, x_k]=\frac a2(\sig_{jk} -s_{jk}), \quad (j\ne k),
$$
$$
[D_j, x_j]=1 + (\iota-1)s_j+\frac a2\sum_{i\ne j}(\sig_{ij}+s_{ij}),
$$
and
$$ \sig_j D_j=-D_j \sig_j, \quad s_{ij}D_{j}=D_i s_{ij}, \quad
\sig_{ij} D_j=-D_i \sig_{ij}.
$$
Consequently
\begin{equation*}
\begin{split}
&\quad (\prod_{j=1}^rD_j^2)(\prod_{j=1}^rx_j^2)\\
&=
(\prod_{j=1}^rD_j)(\prod_{j=1}^rD_j)
(\prod_{j=1}^rx_j)(\prod_{j=1}^rx_j)\\
&=(\prod_{j=1}^rD_j)
\prod_{j=1}^r(D_jx_j -\frac a2 \sum_{i<j}(\sig_{ij} +s_{ij}))
(\prod_{j=1}^rx_j)\\
&=\prod_{j=1}^r(D_jx_j -\frac a2 \sum_{i<j}(\sig_{ij} +s_{ij})+1-(\iota-1)\sig_j)
(\prod_{j=1}^rD_j)(\prod_{j=1}^rx_j)\\
&=\prod_{j=1}^r(D_jx_j -\frac a2 \sum_{i<j}(\sig_{ij} +s_{ij})+1-(\iota-1)\sig_j)
\prod_{j=1}^r(D_jx_j -\frac a2 \sum_{i<j}(\sig_{ij} +s_{ij}))  
\end{split}
\end{equation*}
\end{proof}

\begin{lemm+}\label{res-cay}The operator
$\Res(\Del(\partial_x)\Del(x))$, when acting on $W$-invariant
polynomials of the form $f(x_1^2, \dots, x_r^2)$
with $f(y_1, \dots, y_r)$ being a symmetric polynomials 
in $r$ variables, is given in terms of the coordinates
$y_j=x_j^2$ by 
$$\Res(\Del(\partial_x)\Del(x))f(x)
=\prod_{j=1}^r(U_j^{A}+\frac 12(\iota-1)-\frac 12)
\prod_{j=1}^r U_j^{A}
$$
where $U_j^{A}$ is the Cherednik operator
(\ref{chopa}) of type A acting on
functions of $y_1, \dots, y_r$.
\end{lemm+}
\begin{proof} We use the previous lemma
and consider the operator $D_j x_j-\frac a2\sum_{i<j}(\sig_{ij}+s_{ij})$
acting on the functions of the form
 $f(y_1, \dots, y_r)=f(x_1^2, \dots, x_r^2)$
; this operator maps $f$ into functions of the same form, and
in terms of the variables $y_j=x_j^2$, it
is
\begin{equation}
 \label{dj-c-a}
2(D_j^{(y)} -\frac a2\sum_{i<j}s_{ij} +\frac{\iota-1}2 -\frac 12)
=2(U_j^A+\frac{\iota-1}2 -\frac 12)  
\end{equation}
by direct calculation, proving our Lemma.
\end{proof}

The above lemma
can  also be proved by using
 \cite{Dunkl-cmp-98}, Proposition 5.2,
where the operator $(\prod_{j=1}^r D_j)(\prod_{j=1}^r x_j)$
is expressed in terms of the Cherednik
operators. (Note that our operator $U_j^A$ differs by a constant
with the operator $U_{A, j}$ there.)
Together with (\ref{jack-def}) it  implies then

\begin{lemm+} The Jack symmetric polynomials
 $\Ome_{\m} (x_1^2, \dots, x_r^2)$ are eigenfunction of the operator 
$\Res(\Del^{-\a} \Del(\partial)\Del^{1+\a})$
with eigenvalue
\begin{equation}
  \label{eig-val-1}
\prod_{j=1}^r(\frac a2(r-j) +\frac{\iota}2 +m_j +\a)
(\frac a2(r-j) +m_j +\a).  
\end{equation}
\end{lemm+}

\begin{prop+} For each $\m=(m_1, \dots, m_r)$
there exists a unique polynomial $p_{\n}$
in the space $\mathcal P_{\n}^L$  with $\n$ given by $\m$ as in 
(\ref{n-m-rel-c}) such that
\begin{equation}
  \label{def-pm-c}
  \Res p_{\n}(x_1, \dots, x_r)=\Ome_{\m}(x_1^2, \dots, x_r^2).
\end{equation}
\end{prop+}

\begin{proof} The existence of such polynomial
 $p_{\n}$
in $\PV^L$ is by the Chevalley restriction
theorem. Lemma 4.3 implies that  $p_{\n}$
is in the space $\Pn$, with $\n$ as
given, if and only if
it is an eigenfunction of the operator
$\Del^{-\a}\Del(\partial)\Del^{1+\a}$ with
eigenvalue
\begin{equation*}
 \label{eig-val-2}
\begin{split}
&\qquad \prod_{j=1}^{2r}(\frac{a^\prime}2(2r-j) +1 +n_j+\a)\\
&=\prod_{j=2i-1, i=1}^{r}(\frac{a^\prime}2(2r-2i+1) +1 +n_{2i-1}+\a)
\prod_{j=2i, i=1}^{r}(\frac{a^\prime}2(2r-2i) +1 +n_{2i}+\a)\\
&=\prod_{i=1}^{r}({a^\prime}(r-i+\frac 12) +1 +m_{i}+\a)
\prod_{i=1}^{r}({a^\prime}(r-i) +1 +m_{i}+\a).
\end{split}\end{equation*}
Using 
(\ref{r-r-C}) we see that the factors in the product
  (\ref{eig-val-1}) 
are
$$
\frac a2(r-j)+\frac{\iota}2 +m_j +\a
=a^\prime(r-j)+ 1+\frac{a^\prime}2  +m_j +\a
=a^\prime(r-j+\frac 12)+ 1 +m_j +\a
$$
and
$$
\frac a2(r-j)+ 1 +m_j +\a
=a^\prime(r-j)+ 1+ m_j +\a,
$$
so that   (\ref{eig-val-1}) coincides
with the previous formula.
\end{proof}

We calculate now the Fock space  norm of $p_{\n}$ by
a direct calculation using certain recurrence formula
of Macdonald. 

\begin{lemm+}In terms of the coordinates $y_j=x_j^2,
j=1, \dots, r$, the operator $F_0$ 
when acting on $W$-invariant
polynomials has the following form, 
\begin{equation}
  \label{F0-F0a}
F_0=-(-2F_0^{A} + (\frac{a}2(r-1) +\frac 12 +\frac{\iota-1}2)
\sum_{j=1}^r\partial_j^{(y)})
\end{equation}
where $F_0^A$ is the corresponding operator
for the root system of type A,
\begin{equation}
  \label{F0a}
F_0^{A}=-\frac 12(\sum_{j=1}^r (D_j^A)^2)
=-\frac 12\left(\sum_{j=1}^r y_j (\partial_j^{(y)})^2 +
a\sum_{i\ne j}
\frac{y_i}{y_i-y_j}\partial_i^{(y)}\right)
\end{equation}
\end{lemm+}

With some abuse of notation we denote $\Ome_{\m}^{(x)}$ the function
$\Ome_{\m}^{(x)}=\Ome_{\m}(x_1^2, \dots, x_r^2)$.

To state our next result
we let $\binom{\m}{\m^\prime}$ be the generalized binomial
coefficients \cite{Lassalle-binom} and let
${\m_{j}}$ respectively ${\m^{j}}$
 stand for the
signature ${\m_{j}}=(m_1, \cdots, m_r)-(0, \dots, 0, 1, 0, \dots, 0)$
and ${\m^{j}}=(m_1, \cdots, m_r)+(0, \dots, 0, 1, 0, \dots, 0)$
 (with $1$
in the $j$'th position). The binomial
coefficients  $\binom{\m}{\,\,\m^\prime}$ is then
(loc. cit., Section 14)
\begin{equation}
  \label{binom-co}
\binom{\m}{{\,\,\m_{j}}}  
=(m_j +\frac a2(r-j))\prod_{i\ne j}\frac{m_j-m_i +\frac a2(i-j-1)}
{m_j-m_i +\frac a2(i-j)}
\end{equation}

\begin{lemm+} The operators $E_0$ and $F_0$ have following
 upper  respectively lower triangular form when acting on the
polynomials $\Ome_{\m}^{(x)}$
\begin{equation} 
 E_0\Ome_{\m}^{(x)}
=\sum_{j=1}^r c_{\m}(j)\Ome_{\m^j}^{(x)},
\end{equation}
\begin{equation}
  F_0\Ome_{\m}^{(x)}
=-\frac 14\sum_{j=1}^r \left(4(m_j-1 -\frac a2(j-1)) +
2a(r-1)+ 2(\iota-1) +2 \right)\binom{\m}{\,\,\m_{j}}\Ome_{\m_j}^{(x)},
\end{equation}
where
$$
c_{\m}(j)=\prod_{i\ne j}\frac{m_j-m_i +\frac a2(1+i-j)}
{m_j-m_i +\frac a2(i-j)}.
$$
\end{lemm+}
\begin{proof}This follows from our formula
  (\ref{F0-F0a}) and the result of Macdonold \cite{Mac-leidennotes}, Section D. Following temporarily the notation there, let
$$
\square_1=\sum_{j=1}^r y_j (\partial_j^{(y)})^2 +
a\sum_{i\ne j}
\frac{y_i}{y_i-y_j}\partial_i^{(y)},
$$
and 
$$
\varepsi_1=\sum_{j=1}^r \partial_j^{(y)};
$$
it is proved by Macdonald that
$$
\square_1 \Ome_{\m}
=\sum_{j=1}^r(m_j-1-\frac a2(j-1))
\binom{\m}{{\m_{j}}}
\Ome_{\m_{j}},
$$
$$
\varepsi_1 \Ome_{\m}
=\sum_{j=1}^r\binom{\m}{{\m_{j}}}
\Ome_{\m_{j}}
$$
and
$$
(\sum_{j=1}^r y_j) \Ome_{\m}
=\sum_{j=1}^r c_{\m}(j)
\Ome_{\m_{j}},
$$
which then imply our result.
\end{proof}

\begin{theo+} 
\label{norm-pn-c}
 The   norm square $\Vert p_{\n}\Vert_{\mathcal F}^2$
of  the $L$-invariant polynomial $p_{\n}$,
with $\n$ determined by $\m$ as in (\ref{n-m-rel-c}), is given
by
\begin{equation*}
\begin{split}
&\quad \prod_{1\le i<j\le r}
\frac{\Gamma(\frac a2(j+1-i))
(\frac a2(j-i))}{\Gamma(1+\frac a2(j-1-i)) }
\prod_{j=1}^r(\frac{\iota-1}2+\frac 12 +\frac a2(r-j))_{m_j}
\prod_{j=1}^r(1+\frac a2(r-j))_{m_j}\\
&\times \prod_{1\le i<j\le r}\frac{\Gamma(m_i+1-m_j+\frac a2(j-1-i)) }
{\Gamma(m_i-m_j+\frac a2(j+1-i))
(m_i-m_j+\frac a2(j-i))}  
\end{split}
\end{equation*}
\end{theo+}

\begin{proof}
For a fixed $j$ 
  write $\m^\prime=\m_j=(m_1,\dots,  m_{j-1}, m_j -1, m_{j+1},
\dots, m_r)$ and let $\n^\prime$ be the corresponding
$\n$.  $ p_{\n}$ and  $p_{\n^\prime}$ are orthogonal
in $\mathcal F$ as they are in the different
Schmid components, so are the polynomials
 $ \Ome_{\m}^{(x)}$ and  $\Ome_{\m^\prime}^{(x)}$ 
with the $\Sig$-inner product
by Proposition 5.2. Now, on the
space of all $W$-invariant
polynomials, the adjoint $(E_0)^\ast=-F_0$ with respect
to the $\Sig$-inner product, from this we obtain
$$
(E_0 \Ome_{\m^\prime}^{(x)},  \Ome_{\m}^{(x)})_{\Sig}
-(\Ome_{\m^\prime}^{(x)}, F_0\Ome_{\m}^{(x)})_{\Sig}, 
$$
this gives
$$
\frac{\Vert\Ome_{\m}^{(x)} \Vert_{\Sig}^2}
{\Vert\Ome_{\m^\prime}^{(x)} \Vert_{\Sig}^2}
=\binom{\m}{\m^\prime}\frac{m_j-\frac 12 +\frac{\iota-1}
2  +\frac a2(r-j)}{
c_{\m^\prime}(j)}.
$$
This recursion formula,
together with the fact that
$(p_0, p_0)_{\mathcal F}=
(\Ome_0, \Ome_0)_{\Sig}=1$ uniquely
determine the norm. Carrying out
the calculation gives our result.
\end{proof}

\section{$L$-invariant holomorphic polynomials: 
Type $D_r$ ($r\ge 3$)}

The Weyl group $W$ in this case consists
of the signed permutation of vectors
$(x_1, \dots, x_r)$ keeping the product
$x_1\dots x_r$ invariant.
Thus any $W$-invariant polynomial
is of the form $(x_1\dots x_r)^m f(x_1^2, \dots, x_r^2)$,
where $f$ is a symmetric polynomial in $r$ variables.
The Dunkl operators
 are
$$
D_j^D=\partial_j +
\frac a2\sum_{i\ne j}(\frac 1{x_j -x_i}(1-s_{ij})
+\frac1{x_j -x_i}(1-\sig_{ij})).
$$
We are interested in the operator $(\prod_{j=1}^rD_j)(\prod_{j=1}^rx_j)$. Similar to the proof of Lemma 6.1
we have
\begin{equation}
D_r(\prod_{i=1}^r x_i)=(\prod_{i=1}^{r-1} x_i)
(D_r x_r -\frac a2 \sum_{i<r}(s_{ir}+\sig_{ir})),
\end{equation}
and  generally
\begin{equation}
D_j(\prod_{i\le j}x_i)=
(\prod_{i\le j-1}x_i)(D_j x_j -\frac a2 \sum_{i<j}(s_{ij}+\sig_{ij}).
\end{equation}
We define therefore  Cherednik operator $D_j$, $j=1, \dots, r$
for type $D$,
$$
U_j^D= D_jx_j -\frac a2\sum_{i<j}(s_{ij}+\sig_{ij}).
$$
Thus we have
\begin{equation}
  \label{cay-typ-d}
(\prod_{j=1}^r D_j) (\prod_{j=1}^rx_j)=\prod_{j=1}^r U_j
\end{equation}

Moreover, it is easy to prove that
\begin{equation}
  \label{u-p}
U_j^D (\prod_{j=1}^r x_j)
=(\prod_{j=1}^r x_j) (U_j^D+1),
\end{equation}
Therefore we have
\begin{lemm+} We have the
following formula
$$(\prod_{j=1}^r D_j)^m(\prod_{j=1}^r x_j)^m
=\prod_{j=1}^r (U_j^{D})_m
$$
where
$(T)_m=T(T+1)\dots (T+m-1)$ is the Pochhammer symbol
of any operator $T$.
\end{lemm+}

Consider the Cherednik operator $U_j^D$ acting on
even functions $f(x_1^2, \dots, x_r^2)$. Performing
the change of variables $y_j=x_j^2$, we have
\begin{equation}
U_j^D =2(U_j^A-\frac 12)
\end{equation}
where $U_j^A$ is the operator $(6.2)$. 
This proves consequently  that the Cherednik operators
$U_j^D$ are also commuting operators when acting on even functions.

The next result follows then immediately from
Corollary 5.4, the formulas (7.3) and (7.4).
\begin{lemm+}\label{res-cay-d}The operator
$\Res(\Del(\partial_x)\Del(x))$, when acting on even 
polynomials of the form $f(x_1^2, \dots, x_r^2)$, is,
after the change of variables $y_j=x_j^2$,
given by
$$\Res(\Del(x)^{-\a}\Del(\partial_x)\Del(x)^{1+\a})
=2^{r}\prod_{j=1}^r(U_j^A-\frac 12 +\frac{\a}2)
$$
where
$U^A_j=D_j^{(y)}y_j -\frac a2 \sum_{i<j}\sig_{ij} 
$
is the Cherednik operator in variables $y_j$'s.
\end{lemm+}

\begin{prop+}
 The   polynomials 
$(\prod_{j=1}^rx_j)^m
\Ome_{\m} (x_1^2, \dots, x_r^2)$ are eigenfunction of the operator 
$\Res(\Del(x)^{-\a} \Del(\partial_x)\Del(x)^{1+\a})$
with eigenvalue
$$
\prod_{j=1}^r(a(r-k)+2 +{\a} +
{m} + 2m_j).
$$
Moreover, it is the restriction of a unique
$L$-invariant
polynomials $p_{\n}$ in the space $\Pn$, 
\begin{equation}
  \label{def-pm-d}
  \Res p_{\n}(x_1, \dots, x_r)=(\prod_{j=1}^rx_j)^m
\Ome_{\m}(x_1^2, \dots, x_r^2),
\end{equation}
with $\n$ given by $\m$ as in
(\ref{n-m-rel-d}).
\end{prop+}
\begin{proof} The operator
$\Res(\Del(x)^{-\a} \Del(\partial_x)\Del(x)^{1+\a})$
has the form
$$
(\prod_{j=1}^rx_j)^{-\a}
(\prod_{j=1}^r D_j)(\prod_{j=1}^rx_j)^{1+\a}
$$
by  Corollary 5.4. We calculate its action on 
$(\prod_{j=1}^rx_j)^m
\Ome_{\m}(x_1^2, \dots, x_r^2)$. It is
$$
(\prod_{j=1}^rx_j)^{m}\left((\prod_{j=1}^rx_j)^{-\a-m}
(\prod_{j=1}^r D_j)(\prod_{j=1}^rx_j)^{1+\a+m}\right)
\Ome_{\m}(x_1^2, \dots, x_r^2).
$$
Applying Lemma 7.2 for the operator in the parenthesises 
and using (6.4) we see
that our polynomials is indeed an eigenfunction
with eigenvalue
$$
2^r\prod_{j=1}^r(\frac a2(r-k)+\frac 12 +\frac{\a}2 +
\frac{m}2 + m_j)
$$
as claimed. The rest of the proof is
similar to that of Proposition 6.6.
\end{proof}

\begin{theo+} \label{norm-pn-d}The norm square
$\Vert p_{\n}\Vert_{\mathcal F}^2$
of the $L$-invariant polynomial $p_{\n}$ is given
by
\begin{equation*}
  \label{norm-pn-dd}
\begin{split}
&\quad
 \prod_{1\le i<j\le r}\frac{\Gamma(\frac a2(j+1-i))
(\frac a2(j-i))}
{\Gamma(1+\frac a2(j-1-i)) }
 2^{2(m_1+\dots +m_r)}
\prod_{j=1}^r(\frac 12 +\frac a2(r-j))_{m_j}
\prod_{j=1}^r(1+\frac a2(r-j))_{m_j}
\\
&\quad \times
\prod_{j=1}^r (\frac a2(r-j)+1+m_j -\frac 12)_m  
\prod_{1\le i<j\le r}\frac{\Gamma(m_i+1-m_j+\frac a2(j-1-i)) }
{\Gamma(m_i-m_j+\frac a2(j+1-i))
(m_i-m_j+\frac a2(j-i))}
\end{split}
\end{equation*}
\end{theo+}
\begin{proof} Proposition 5.2 implies that
$$\Vert p_{\n}\Vert_{\mathcal F}^2=
\Vert\Res p_{\n}\Vert_{\Sig}^2=
\Ome_{\m}(D_1^2, \dots, D_r^2) (\prod_{j=1}^r D_j)^m(\prod_{j=1}^r x_j)^m
\Ome_{\m}(x_1^2, \dots, x_r^2)\big{\vert}_{x=0}.
$$
However, by Lemma 7.1
$$(\prod_{j=1}^r D_j)^m(\prod_{j=1}^r x_j)^m
=\prod_{j=1}^r (U_j)_m=2^{r}\prod_{j=1}^r (U_j^A -\frac 12)_m,
$$
which has $\Ome_{\m}(x_1^2, \dots, x_r^2)$
as an eigenfunction with eigenvalue
$$
\prod_{j=1}^r (\frac a2(r-j)+1+m_j -\frac 12)_m
$$
by  (\ref{jack-def}).
Thus
\begin{equation*}
\begin{split}
\Vert p_{\n}\Vert_{\mathcal F}^2&=
\prod_{j=1}^r (\frac a2(r-j)+1+m_j -\frac 12)_m
\Ome_{\m}(D_1^2, \dots, D_r^2) \Ome_{\m}(x_1^2, \dots, x_r^2)\big{\vert}_{x=0}\\
&=\prod_{j=1}^r (\frac a2(r-j)+1+m_j -\frac 12)_m
(\Ome_{\m}^{x}, \Ome_{\m}^{x})_{\Sig}.  
\end{split}
\end{equation*}
The inner product above can be calculated by
the same recursion formula as for type C. Indeed
the Laplace operator $F_0$ in this case, after changing of variables
is of the form 
$$
F_0=-2\left(-2F_0^{A} +
(\frac{a}2(r-1) +\frac 12 )\sum_{j=1}^r\partial_j^{(y)}\right)
$$
which is of the same
form as   (\ref{F0-F0a}) except the term
$\frac{\iota-1}2$ is missing and that the coefficient
$-1$ in front is replace by $-2$.
\end{proof}

\section{Bargmann transform of $L$-invariant
polynomials: Flat case}

Consider again the restriction map 
\begin{equation}
R=R_{\nu}: \mathcal F_\nu (V_{\bc}) \mapsto C^{\omega}(V),
\quad Rf(x)=f(x)e^{-\frac{\nu}2\Vert x\Vert^2 }
\end{equation}
It defines a bounded operator and the isometric
part $U$ in the polar decomposition $R=|R|U$ of
$R$ is a unitary operator, and is the Segal-Bargmann transform.
The Berezin transform in this case
is
$$
|R|^2f(x)= 
RR^\ast f(x)=\int_V e^{-\frac \nu 2 \Vert x-y\Vert^2 }f(y) dy
$$

Recall our identification of $V$ with $\fq$. In this
identification  $\fa=\br  e_1 +\dots +\br e_r$
is a subspace of $V$.

Define the Bessel function on $V$ by
\begin{equation}
  \label{Bes-func}
J_{\blam}(x)=\int_{L} e^{-i\blam((lx))}  dl
\end{equation}
for $\blam \in \fa^\ast$ where
we extend the linear functional $\blam$ on
$\fa$ to $V$ via the orthogonal
projection onto $\fa$ and where $dl$ is the normalized Haar measure
on $L$, so that $J_{\blam}(0)=1$.

\begin{defi+} The Hermite polynomial $\zeta_{\n, \nu}(\blam)$
related to the root system $\Sig$ of type C or D on $\fa$ is defined by
 the Rodrigues type formula:
$$
\zeta_{\n, \nu}(\blam)=
\frac{1}{\Vert p_{\n}\Vert_{\mathcal F}^2}
p_{\n}(\partial_x)
\left(e^{\frac{\nu}2 \Vert x\Vert^2}J_{\blam}(x)\right)\big{\vert}_{x=0}.
$$
for those $\n$ 
as determined in (\ref{n-m-rel-d}). Here $p_{\n}(\partial_x)$
is the differential operator on $V$ obtained from
$p_{\n}(x)$ by the same convention
as in   (\ref{fock-nu}).
\end{defi+}
Conceptually it is better to write this
as
$$
\zeta_{\n, \nu}(\blam)=\frac{1}{\Vert p_{\n}\Vert_{\mathcal F}^2}
p_{\n}(\partial_z)(R^{-1}J_{\blam})(0)
=\frac{1}{\Vert p_{\n}\Vert_{\mathcal F}^2}
p_{\n}(\partial_z)(e^{\frac{\nu}2(z, \bar z)}
J_{\blam}(z))(0)
$$
where we extend the function $e^{\frac{\nu}2\Vert x\Vert^2}$ on $V$
to a \textit{holomorphic} function $e^{\frac{\nu}2(z, \bar z)}$
 on the whole
space $V_{\bc}$.

The following expansion is a direct consequence 
of the definition.

\begin{lemm+}
 The $L$-invariant
function $e^{\frac{\nu}2\Vert x\Vert^2}J_{\blam}(x)
$ on $V$ has the following expansion in terms of $p_{\n}(x)$
\begin{equation}
  \label{exp-bes}
  e^{\frac{\nu}2\Vert x\Vert^2 }J_{\blam}(x)
=\sum_{\n} p_{\n}(x) \zeta_{\n, \nu}(\blam),
\end{equation}
where the summation is over all $\n$ as determined
in (\ref{n-m-rel-c}) and (\ref{n-m-rel-d}).
\end{lemm+}

The next result computes the Fourier
transform of the  Segal-Bargmann transform $Up_{\n, \nu}$ of 
$ p_{\n}(\blam)$.
Denote $f(x)\mapsto \widetilde{f}(\blam)$ the
Fourier transform on $V$ evaluated on $\fa^\ast$,
\begin{equation}
  \label{fouri-def}
  \widetilde{f}(\blam)=\int_V f(x)e^{i\blam(x)}dx, \qquad \blam\in  
\fa^\ast.
\end{equation}

\begin{prop+} The Fourier transform 
$ \widetilde{Up_{\n}}(\blam) $ of
$Up_{\n}$ is
\begin{equation}
  \label{xi-up-f}
  \widetilde{Up_{\n}}(\blam)  =(-\nu)^{-|\n|}
(\frac {2\pi}{\nu})^{\frac d4} \zeta_{\n, \nu}(\blam) e^{-\frac 1{4\nu}
\Vert \blam\Vert^2}
\end{equation}
\end{prop+}

\begin{proof}The Berezin transform is
the convolution operator with Gaussian kernel,
whose Fourier transform is also a Gaussian. Namely
$|R|^2$ has $e^{i\blam(x)}$ and therefore the Bessel function as its eigenfunction. More precisely, we have
\begin{equation}
  \label{gau-fou}
|R|^2J_{\blam}(x)=
(\frac {2\pi}{\nu})^{\frac{d}2} e^{-\frac 1{2\nu}
\Vert \blam\Vert^2} J_{\blam}(x).
\end{equation}
We may rewrite it  as
$$
\int_{V}e^{\nu (x, y)}e^{-\frac \nu 2\Vert y\Vert^2}J_{\blam}(x)
dy
=(\frac {2\pi}{\nu})^{\frac{d}2}
 (J_{\blam}(x) e^{\frac{\nu} 2\Vert x\Vert^2})
e^{-\frac 1{2\nu}
\Vert \blam\Vert^2}.
$$
We differentiate both side by the differential
operator $p_{\n}(\partial_x)$ and evaluate
at $x=0$. To do this, we observe that
$$
p_{\n}(\partial_x)
e^{-\nu (x, y)}(0)
=(-\nu)^{|\n|}p_{\n}(y),
$$
so that the resulting formula is
\begin{equation}
  \label{tric}
(-\nu)^{|\n|} \int_{V}p_{\n}(y)e^{-\frac \nu 2\Vert y\Vert^2}J_{\blam}(y)
dy=
(\frac {2\pi}{\nu})^{\frac{d}2} 
\Vert p_{\n}\Vert_{\mathcal F}^2
\zeta_{\n, \nu}(\blam) e^{-\frac 1{2\nu}
\Vert \blam\Vert^2}.
\end{equation}
The left  hand side
is actually
\begin{equation}\label{fou-rf}
(-\nu)^{|\n|}\widetilde{Rp_{\n}}(\blam)
=(-\nu)^{|\n|}\widetilde{|R|U p_{\n}}(\blam).
\end{equation}
On the other hand formula 
the formula (\ref{gau-fou}) implies that
 $$
\widetilde{|R|^2f}(\blam)=
(\frac {2\pi}{\nu})^{\frac d 2} e^{-\frac 1{2\nu}
\Vert \blam\Vert^2} \widetilde{f}(\blam)
$$
for any $f\in L^2(V)$.
 Thus
$$
\widetilde{|R|f}(\blam)=
(\frac {2\pi}{\nu})^{\frac d4}  e^{-\frac 1{4\nu}
\Vert \blam\Vert^2}
 \widetilde{f}(\blam).
$$
Substituting this into (\ref{fou-rf}) for 
$f=Up_{\n}$ we get
$$
(-\nu)^{|\n|}\widetilde{Rp_{\n}}(\blam)
=(-\nu)^{|\n|}(\frac {2\pi}{\nu})^{\frac d4}  e^{-\frac 1{4\nu}
\Vert \blam\Vert^2}\widetilde{Up_{\n}}(\blam).
$$
The equality (\ref{tric}) becomes now
$$
(-\nu)^{|\n|}(\frac {2\pi}{\nu})^{\frac d4}  e^{-\frac 1{4\nu}
\Vert \blam\Vert^2}\widetilde{Up_{\n}}(\blam)  
=(\frac{2\pi}{\nu})^{\frac d2}\Vert p_{\n}\Vert^2_{\mathcal F}
 \zeta_{\n, \nu}(\blam) e^{-\frac 1{2\nu}
\Vert \blam\Vert^2},
$$
namely
\begin{equation} \label{xi-up-b}
  \widetilde{Up_{\n}}(\blam)  =
(\frac{2\pi}{\nu})^{\frac d4}
(-\nu)^{-|\n|}\Vert p_{\n}\Vert^2_{\mathcal F}
 \zeta_{\n, \nu}(\blam) 
e^{-\frac 1{4\nu}\Vert \blam\Vert^2},
\end{equation}
completing the proof.
\end{proof}

The unitarity of $U$ and the Fourier transform
implies then  the orthogonality relation
of $\zeta_{\n, \nu}$. Let  $C_1$
be the normalization constant 
so that the following
measure is a probability measure on $\fa$:
\begin{equation}
\label{c1-c}
C_1 (\frac{2\pi}{\nu})^{\frac d2}
e^{-\frac 1{2\nu}\Vert\blam\Vert^2}\prod_{j=1}^r{|\lam_j|^{\iota-1}}
\prod_{1\le i<j\le r}^r{|\lam_i-\lam_j|^{a}}
\end{equation}
for root system of type $C$,
and \begin{equation}
\label{c1-d}
C_1 (\frac{2\pi}{\nu})^{\frac d2}
e^{-\frac 1{2\nu}\Vert\blam\Vert^2}\prod_{1\le i<j\le r}^r{|\lam_i^2-\lam_j^2|^{a}}
\end{equation}
for type $D_r$.
 (Note $C_1$  is independent
to the parameter $\nu$.) $C_1$ 
can be evaluated by the Selberg-Macdonald
formula (proved by Opdam \cite{Opdam-bessel}). 

\begin{coro+} The $W$-invariant polynomials
$\Vert p_{\n}\Vert_{\mathcal F_{\nu}}
{\zeta_{\n, \nu}(\blam)}$ for an orthonormal
basis for the Hilbert space of $W$-invariant
$L^2$-functions with the above
probability measure.
\end{coro+}

The expansion  (\ref{exp-bes}) takes now 
the following form
\begin{coro+} The $L$-invariant
function $e^{\frac{\nu}2\Vert x\Vert^2}J_{\blam}(x)
$ on $V$ has the following expansion in terms of $p_{\n}(x)$,
writing $q_{\n, \nu}=
\frac{p_{\n}(x)}{\Vert p_{\n}\Vert_{\mathcal F_\nu}}$
the normalized $L$-invariant polynomial,
\begin{equation}
  \label{exp-bes-new}
e^{-\frac 1{4\nu}
\Vert \blam\Vert^2}  e^{\frac{\nu}2\Vert x\Vert^2}J_{\blam}(x)
=(\frac {2\pi}{\nu})^{\frac d2} 
\sum_{\n} (-1)^{|\n|}{q_{\n, \nu}(x)}\widetilde {U q_{\n, \nu}}
(\blam).
\end{equation}
\end{coro+}

\section{Segal-Bargmann transform of $L$-invariant
polynomials: bounded case}

The result is this section is parallel to
 that of the previous section
and is an explicit realization of the
Corollary \ref{realize!}. We will be brief;
see also \cite{jp+gkz-pl3} and \cite{oz-weyl}. 

Recall the Berezin transform $B_{\nu}$ in (\ref{def-Bnu}).
It is proved in \cite{gkz-bere-rbsd} that $B_{\nu}$ 
defines a $H$-invariant bounded operator
on $L^2(D)$ with the invariant measure (\ref{H-invmea}). Let
$b_{\nu}(\blam)$ be its spectral symbol, namely
$$
B_{\nu}\phi_{\blam} = b_{\nu}(\blam)\phi_{\blam}
$$
in the sense of spectral decomposition, where
$\phi_{\blam}$ on $D$ is the spherical function;
the function $b_{\nu}(\blam)$ is explicitly
calculated in \cite{gkz-bere-rbsd}.

Define the polynomials $\xi_{\n, \nu}(\blam)$
by the Rodrigues type formula
$$
\xi_{\n, \nu}(\blam)=
\frac{1}{\Vert p_{\n}\Vert_{\mathcal F}^2}
p_{\n}(\partial_x)(h^{-\frac \nu 2}(x)\phi_{\blam}(x))\big{\vert}_{x=0}.
$$
Thus $\xi_{\n, \nu}(\blam)$ is $W$-invariant. 

\begin {theo+}  \label{exp}
The $L$-invariant analytic function $h(x, \bar x)^{-\frac \nu 2}\phi_{\blam}(x)$,
when extended to a holomorphic function 
$h(z, z)^{-\frac \nu 2}\phi_{\blam}(z)$ 
in
a neighborhood of $D$ in $\mathcal D$,
has the following
expansion
in terms of the $L$-invariant
 polynomials $p_{\n}(z)$,
\begin{equation}
\label{exp-sph}
h(z, z)^{-\frac \nu 2}\phi_{\blam}(z)
=\sum_{\n} \xi_{\n, \nu}(\blam){p_{\n}(z)}
\end{equation}
near $z=0$ in $\mathbb D$. Moreover the spherical transform
of the Segal-Bargmann transform of $p_{\n}$ is
\begin{equation}
\label{barg-pm}
\widehat{Up_{\n}}(\blam)=\frac1{c_{\nu}^{\frac 12}}b_{\nu}(\blam)^{\frac 12}
\xi_{\n, \nu}(\blam)\Vert p_{\n}\Vert_{\nu}^2.
\end{equation}
Thus $\Vert p_{\n}\Vert_{\nu}\xi_{\n}(\blam)$ for all $\n$ form an
orthonormal basis for the
space $L^2(\fa^\ast, C_0 \frac{b_{\nu}(\blam)}{c_{\nu}}|c(\blam)|^{-2})^W$.
\end{theo+}

Conceptually it is better to write
(\ref{exp-sph}) in the form
\begin{equation}
\label{exp-sph-new}
(R^{-1}\phi_{\blam})(z)=h(z, z)^{-\frac{\nu}2}\phi_{\blam}(z)
=\sum_{\n}({\Vert p_{\n}\Vert_{\nu}}
{ \xi_{\n, \nu}(\blam)})
\otimes \frac{p_{\n}(z)}{\Vert p_{\n}\Vert_{\nu}}.
\end{equation}
Note that we have extended the real analytic function
$h(x, \bar x)^{-\frac{\nu}2}$ to  a \textit{holomorphic}
function, formally written as
 $h(z, z)^{-\frac{\nu}2}$, on $\mathbb D$. 
Notice also that both  $\frac{p_{\n}(z)}{\Vert p_{\n}\Vert_{\nu}}$
and ${\Vert p_{\n}\Vert_{\nu}}
{ \xi_{\n}(\blam)}$ are orthonormal basis
in the respective Hilbert spaces. Thus $(R^{-1}\phi_{\blam})(z)$
is the Schwarz kernel
for the unitary operator 
from $\mathcal H_{\nu}^L$
onto $L^2(\fa^\ast, C_0 \frac{b_{\nu}(\blam)}{c_{\nu}}|c(\blam)|^{-2})^W$
obtained by taking the composition 
 of the spherical and Segal-Bargmann transforms.

\begin{rema+} Associated to each root
system there are the multi-variable Askey-Wilson polynomials
determined by the root multiplicities (which
are positive real numbers)
and two more extra parameters \cite{vanDiejen-tams}. Our
 polynomial $\xi_{\n, \nu}(\blam)$
depend on  extra parameter
and can be reviewed as some limiting cases of
the Askey-Wilson polynomials polynomials.  We have
thus found an alternative
simple proof of the their
orthogonality relation.
Moreover
the (known) expansion of $h^{-\frac \nu 2}$
in terms of the Jack symmetric polynomials
will also give some evaluation formula for
the polynomials; we will however not pursue
it here.
\end{rema+}

\begin{rema+} Consider the spherical function
on the Hermitian  symmetric space $\mathbb D$
or on its compact dual. There are Jacobi type functions.
There have been attempts (see \cite{Beerends-Opdam})
in expanding
the Jacobi  functions in terms of the Jack symmetric
polynomials and study the combinatorial
properties of the coefficients.
 It  turned out (see \cite{gkz-invtoe}, \cite{oz-weyl}
and \cite{gz-invdiff}) however that the expansion
of the spherical transform multiplied with
the reproducing kernel has a much better
analytical significance. This is also the case here.
Further more by considering analytic continuation
in the parameter $\nu$ we may recover
the expansion of the spherical functions itself.
\end{rema+}


\section*{Appendix 1: Evaluation
of the constant $C_0$}

We consider type $D_r$ ($r\ge 3$) first, in this case $\Del(x)=\del(x)$
and the rank $r(\mathbb D)$ is $r$.

It is proved in \cite{Hel2}, Chapter IV, Exercise C4, that when
the $H$-invariant measure  $d\iota(z)$ on
$D=H/K$ is (regularly) normalized so that 
$$
\int_Df(z)d\iota(z)=\int_{\fa^+}
f(\exp H)\prod_{\a\in \Sig^+}(e^{\a(H)}
-e^{-\a(H)})^{m_\a} d_nH
$$
with $d_nH$ the regular normalization,
then the Plancherel formula reads
$$
|W|\int_{D}|f(z)|^2d\iota(z)
=\int_{\fa^\ast}|\widehat f(\blam)||c(\blam)|^{-2}d_n\blam.
$$
By regular  normalizations $d_nH$
 on $\fa$  it is meant that the $d_nH$ is
the Euclidean measure on $\fa$ induced
by the restriction on it of the
Killing form on $\fg$, multiplied with the
factor $\frac{1}{\sqrt {2\pi}^r}$; by the Riesz lemma 
we get an identification
of $\fa^\ast$ with $\fa$ and thus similarly get
a regular measure $d_n\blam$ on $\fa^\ast$.
In our case they
are
$$
d_nH=\frac 1{\sqrt{2\pi}^{r}}   {\sqrt{2a(r-1)}^r}
dx_1\cdots dx_r,
\quad d_n\blam =\frac{1}{\sqrt{2\pi}^{r}}
\left(\frac{2}{\sqrt{2a(r-1)}}\right)^r
d\lam_1\cdots d\lam_r,
$$
if $H=x_1\xi_1+\dots x_r\xi_r\in \fa^\ast$ and 
$\blam =\lam_1\beta_1+\dots
+\lam_r\beta_r\in \fa^\ast$ and the order
$|W|$ of the Weyl group is $2^{r-1}r!$. 
There exists  now a constant $C_0$ so that
$$
\int_D f(z)\frac{dm(z)}{h^{\frac p2}(z)}
=C_0\int_D f(z)d\iota(z)
=C_0\int_{\fa^+}
f(\exp H)\prod_{\a\in \Sig^+}(e^{\a(H)}
-e^{-\a(H)})^{m_\a} d_nH.
$$
Take $f(z)=h(z, \bar z)^{\sig}$ for sufficiently large $\sig$.
The left hand side can, by performing Cayley
transform $z\mapsto w=\frac{e+z}{e-z}$, 
$z=\frac{w-e}{w+e}$ mapping
$D$ to its the Siegel domain $\mathcal S$
(see \cite{Loos-bsd}, \cite{gkz-bere-rbsd}),  be evaluated
eventually  by the Gindikin Gamma function, while the right
hand side is a kind of Selberg type integral and
is also known.  Let us carry out this calculation. 
Write $I_1$ for the left hand side,
$$
I_1=\int_D h^{\sig-\frac p2}(z){dm(z)}
=4^{r(\sig-\frac p2)}  2^{d} \int_{\mathcal S}
\frac{\Del(w)^{\sig-\frac p2}}
{\Del(e+w)^{2\sig-p}}   \frac{dw}{\Del(w+e)^{\frac dr}}
$$
since the determinant of the
differential of the Cayley transform
$z=\ga_{-1}(w)$ is 
$$J_{\ga_{-1}}(w)=2^d\Del(e+w)^{-2\frac dr},
$$ 
(see \cite{FK-book}, Chapter X, Proposition 2.4,
for the calculation of the complex Jacobian of the Cayley
transform), and that
$$ 
h(z,\bar z)=4^{r}\frac{\Del(w)}
{\Del(e+w)^{2}}.
$$

This integral over
 $\mathcal S$ is evaluated in \cite{gkz-bere-rbsd},
Proposition 4.2, 
 the result is
$$
I_1=4^{r(\sig-\frac p2)}  2^{d}
{\sqrt \pi}^{d} 4^{n_{A}-r(\sig+\frac 12) }  
\frac{\GGa(\nu -\frac{d_B}r)}
{\GGa(\sig +\frac 12)}
={\sqrt \pi}^{d}\frac{\GGa(\sig -\frac{d_B}r)}
{\GGa(\sig +\frac 12)};
$$
here, adopting the notation there, $\GGa(x)$ is
the Gindikin Gamma function.
The right hand side is
$\frac 1{\sqrt{2\pi}^{r}}   {\sqrt{2a(r-1)}^r}
I_0
$
where
$$I_0=\int_{x_1>x_2>\dots > x_r>0}
\prod_{j=1}^r(1-\tanh^2 x_j)^{\sig}
\prod_{i<j}(e^{x_i -x_j}-e^{x_j -x_i})^{a} (e^{x_i +x_j}-e^{-x_i-x_j})^{a} 
dx_1\cdots dx_r.
$$
Perform change of variables $t_j=\tanh x_j$
then 
$dx_j=(1-t_j^2)^{-1}dt_j$
and
$$
(e^{x_i -x_j}-e^{x_j -x_i})(e^{x_i +x_j}-e^{-x_i -x_j})
=\frac{4(t_i^2-t_j^2)}{(1-t_i^2)(1-t_j^2)}.
$$
So that
\begin{equation}
\begin{split}
I_0&=   4^{\frac a2r(r-1)}
\int_{1>t_1>t_2>\dots > t_r>0}
\prod_{j=1}^r(1-t_j^2)^{\sig-a(r-1)-1}
\prod_{i<j}(t_i^2-t_j^2)^a dt_1\cdots dt_r\\
&=  4^{\frac a2r(r-1)} \frac 1{r!}
\int_{[0, 1]^r}
\prod_{j=1}^r(1-t_j^2)^{\sig-a(r-1)-1}
\prod_{i<j}(t_i^2-t_j^2)^a dt_1\cdots dt_r.
\end{split}
\end{equation}
Changing again  variables, letting $s_j=t_j^2$, the integral becomes
\begin{equation}
I_0= 4^{\frac a2r(r-1)}2^{-r} \frac 1{r!} 
\int_{[0, 1]^r}
\prod_{j=1}^r(1-s_j)^{\sig-a(r-1)-1}
\prod_{j=1}^r s_j^{-\frac 12}
\prod_{i<j}(s_i-s_j)^a  ds_1\cdots ds_r;
\end{equation}
its value is
\begin{equation}
\begin{split}
\frac {2^{ar(r-1)-r}} {r!}
\prod_{i=1}^r\frac{\Gamma(\sig-\frac a2(r-1)- \frac a2(i-1))
\Gamma(\frac 12 +\frac a2(r-1)- \frac a2(i-1))
}{\Gamma(\sig +\frac 12 -\frac a2(i-1))}
\frac{\prod_{j=1}^r\Gamma(\frac a2(r-j+1))}
{\Gamma(\frac a2)^r};
\end{split}
\end{equation}
see \cite{Macd-book}, Chapter VI, example 7, pp. 385-386. One of the  fraction
in the above product is the following, and can also be written as
$$
\prod_{i=1}^r\frac{\Gamma(\sig-\frac a2(r-1)- \frac a2(i-1))
}{\Gamma(\sig +\frac 12 -\frac a2(i-1))}
=\frac{\GGa(\sig-\frac a2(r-1))} {\GGa(\sig +\frac 12)}
=\frac{\GGa(\sig-\frac {d_B}r)} {\GGa(\sig +\frac 12)};
$$
this fraction appears also in $I_1$.
 Finally the constant $C_0$ is determined
by the formula
$$
I_1=C_0 \frac{\sqrt{2ar}^r}{\sqrt{2\pi}^{r}} I_0
$$
and we find that
$$
C_0= \frac{\sqrt{2\pi}^{r}} {\sqrt{2ar}^r}
\frac{I_1}{I_0}
=\frac{\sqrt{2\pi}^{r}} {\sqrt{2ar}^r}{\sqrt \pi}^{d}
\frac{r!} {2^{ar(r-1)-r}} 
\frac{\Gamma(\frac a2)^r} {
\prod_{j=1}^r\Gamma(\frac a2(r-j+1))
\Gamma(\frac 12 +\frac a2(r-1)- \frac a2(i-1))
  } 
$$


The constant for type  $C_r$ can be
evaluated similarly and we leave it
to the interested reader.

\section*{Appendix 2: Irreducible real bounded symmetric domains}

We list here the the associated
Lie algebras 
 $(\fg, \fk)$ and $(\fh, \fl)$
 of non-compact irreducible bounded symmetric
domains $\mathbb D=G/K \subset V_{\mathbb C}$ and respectively
their real bounded symmetric subdomain $D=H/L\subset V$,
when $V$ is not an Euclidean Jordan algebra; see \cite{Loos-bsd}. They were also classified by Olafsson
(\cite{Ola-dg}, \cite{Faraut-Olafsson}) through
 Lie theoretic methods.
 The restricted
root systems $\Sig$ of $(\fh, \fl)$ are also indicated,
the name of the types (Types B, BC, etc)
 is in
consistence with Loos \cite{Loos-bsd}, Proposition 11.18.

$\ \underset{}{\text{
Table 1. Irreducible  real bounded symmetric subdomains}}$

\begin{center}
{\begin{tabular}{|l|l|l|l|} 
\hline
 $\mathfrak{g}, \mathfrak{k} $ &$\mathfrak{h}, \mathfrak{l}$ 
&         $\Sig$ \\
\hline 
 $ \mathfrak{su}(r, r), s(\mathfrak{u}(r)\oplus \mathfrak{u}(r))$  &

$ \mathfrak{so}(r, r), \mathfrak{so}(r)\oplus \mathfrak{so}(r)$
 & $D_r$
\\
\hline
 $ \mathfrak{su}(r, r+b), s(\mathfrak{u}(r)\oplus \mathfrak{u}(r+b))$  &

$ \mathfrak{so}(r, r+b), \mathfrak{so}(r)\oplus \mathfrak{so}(r+b)$
 & $B_r $
\\
\hline
 $ \mathfrak{su}(2r, 2r), 
s(\mathfrak{u}(2r)\oplus \mathfrak{u}(2r))$  &

$ \mathfrak{sp}(r, r), \mathfrak{sp}(r)\oplus \mathfrak{sp}(r)$
&$ C_r$
\\
\hline
 $ \mathfrak{su}(2r, 2r+2b), 
s(\mathfrak{u}(2r)\oplus \mathfrak{u}(2r+2b))$  &

$ \mathfrak{sp}(r, r+b), \mathfrak{sp}(r)\oplus \mathfrak{sp}(r+b)$
&$BC_r $
\\
\hline
 $ \mathfrak{so^\ast}(4r), \mathfrak{u}(2r)$  &
$ \mathfrak{so}(2r, \mathbb C), \mathfrak{so}(2r)$
 & $D_{2r}$
\\
\hline
 $ \mathfrak{so^\ast}(2r), \mathfrak{u}(r)$  &
$ \mathfrak{so}(r, \mathbb C), \mathfrak{so}(r)$
 & $B_r \,(r \, \text{odd})$
\\
\hline
 $ \mathfrak{so}(2, p+q), \mathfrak{so}(2)\oplus
\mathfrak{so}(p+q))$  &
$ \mathfrak{so}(1, p)\oplus
\mathfrak{so}(1, q), 
\mathfrak{so}(p)\oplus \mathfrak{so}(q)$
 & $D_2 (q\ne 0) $
\\
\hline
 $ \mathfrak{so}(2, p), \mathfrak{so}(2)\oplus
\mathfrak{so}(p))$  &
$ \mathfrak{so}(1, p), 
\mathfrak{so}(p)$
 & $\quad C_1  $
\\
\hline
 $ \mathfrak{sp}(2r, \mathbb R), \mathfrak{u}(2r)$ 
 &
$ \mathfrak{sp}(r,  \mathbb C),
\mathfrak{sp}(r)$
&$ C_r (b=0)$
\\
\hline
 $ \mathfrak{e}_{6(-14)},  \mathfrak{so}(10)+\mathbb R$ 
 &
$ \mathfrak{sp}(2, 2),\mathfrak{sp}(2)\oplus
\mathfrak{sp}(2)$
&
$ B_2$
\\
\hline
 $ \mathfrak{e}_{6(-14)},  \mathfrak{so}(10)+\mathbb R$ 
 &
$ \mathfrak{f}_{4(-20)},\mathfrak{so}(9)$
&
$ BC_1$
\\
\hline
 $ \mathfrak{e}_{7(-25)},  \mathfrak{e}_6+\mathbb R$ 
 &
$ \mathfrak{su}^\ast(8),\mathfrak{sp}(4)$
&
$ D_3$
\\
\hline

\end{tabular}
\medskip
\medskip
\medskip

}\end{center}

\bibliographystyle{amsplain}

\begin{thebibliography}{10}

\bibitem{Beerends-Opdam}
R.~Beerends and E.~Opdam, \emph{Certain hypergeometric series related to the
  root system {BC}}, Trans. AMS. \textbf{339} (1993), 581--609.

\bibitem{Chered-uni}
I.~Cherednik, \emph{A unification of the {Knizhnik-Zamolodchikov} and {Dunkl}
  operators via affine {Hecke} algebras}, Invent. Math. \textbf{106} (1991),
  411--432.

\bibitem{doz}
M.~Davidson, G.~Olafsson, and G.~Zhang, \emph{{Segal-Bargmann transform on
  Hermitian symmetric spaces}}, in preparation.

\bibitem{Dunkl-tams}
C.~F. Dunkl, \emph{Differential-difference operators associated to reflection
  groups}, Trans. Amer. Math. Soc. \textbf{311} (1989), 167--183.

\bibitem{Dunkl-cjm}
\bysame, \emph{Integral kernels with reflection group invariance}, Canadian J.
  Math. \textbf{43} (1991), 1213--1227.

\bibitem{Dunkl-cmp-98}
\bysame, \emph{Orthogonal polynomials of types{ A} and {B} and related
  {Calogero} models}, Commun. Math. Phys. \textbf{197} (1998), 451--487.

\bibitem{FK}
J.~Faraut and A.~Koranyi, \emph{Function spaces and reproducing kernels on
  bounded symmetric domains}, J. Funct. Anal. \textbf{88} (1990), 64--89.

\bibitem{FK-book}
\bysame, \emph{Analysis on symmetric cones}, Oxford University Press, Oxford,
  1994.

\bibitem{Faraut-Olafsson}
J.~Faraut and G.~Olafsson, \emph{Causal semisimple symmetric spaces, the
  geometry and harmonic analysis}, Semigroups in algebra geometry and analysis
  (Berlin) (Lawson Hofmann and Vinberg, eds.), De Gruyter, 1995.

\bibitem{Heckman-remark}
G.~J. Heckman, \emph{{A remark on the Dunkl differential-difference
  operators}}, Harmonic analysis on reductive groups (Brunswick, ME, 1989)
  (Boston, MA), Progr. Math., vol. 101, {Birkh\"auser Boston}, 1991,
  pp.~181--191.

\bibitem{Hel1}
S.~Helgason, \emph{Differential geometry and symmetric spaces}, Academic Press,
  New York, London, 1978.

\bibitem{Hel2}
\bysame, \emph{Groups and geometric analysis}, Academic Press, New York,
  London, 1984.

\bibitem{Ola-Hilg}
J.~Hilgert and G.~\'Olafsson, \emph{Causal symmetric spaces, geometry and
  harmonic analysis}, Perspectives in Mathematics, vol.~18, Academic Press,
  1996.

\bibitem{Howe-spgp}
R.~Howe, \emph{The first fundamental theorem of invariant theory and spherical
  subgroups}, Algebraic groups and their generalizations: classical methods
  (University Park, PA, 1991), Part I, Proc. Sympos. Pure Math., vol.~56, Amer.
  Math. Soc., Providence, RI, 1994.

\bibitem{HU}
R.~Howe and T.~Umeda, \emph{The {Capelli} identity, the double commutant
  theorem, and multiplicity free action}, Math. Ann. \textbf{290} (1991),
  565--619.

\bibitem{Hua}
L.~K. Hua, \emph{Harmonic analysis of functions of several complex variables in
  the classical domains}, Amer. Math. Soc., Providence, Rhode Island, 1963.

\bibitem{toshi-multfree-rest}
T.~Kobayashi, \emph{Multiplicity-free restrictions of unitary highest weight
  modules for reductive symmetric pairs}, preprint.

\bibitem{toshi-saga}
\bysame, \emph{Multiplicity-free branching laws for unitary highest weight
  modules}, Proceedings of the Symposium on Representation Theory held at
  {Saga}, Kyushu 1997 (K.~Mimachi, ed.), 1997, pp.~8--17.

\bibitem{Kraemer}
M.~{Kr\"amer}, \emph{{Sph\"arische Untergruppen in kompakten
  zusammenh\"angenden Liegruppen}}, Compositio Math. \textbf{38} (1979), no.~2,
  129--153.

\bibitem{Lassalle-binom}
M.~Lassalle, \emph{Coefficients binomiaux g\'en\'eralis\'es et polyn\^omes de
  {Macdonald}}, J. Functional Analysis \textbf{158} (1998), 289--324.

\bibitem{Loos-bsd}
O.~Loos, \emph{Bounded symmetric domains and {Jordan} pairs}, University of
  California, Irvine, 1977.

\bibitem{Mac-leidennotes}
I.~G. Macdonald, \emph{Hypergeometric functions}, Lecture notes, unpublished.

\bibitem{Macd-book}
\bysame, \emph{Symmetric functions and {Hall} polynomials}, second ed.,
  Clarendon Press, Oxford, 1995.

\bibitem{Neretin-beta-int}
Yu. Neretin, \emph{Matrix analogs of the integral {$B(\a, \rho-\a)$} and
  {Plancherel} formula for {Berezin} kernel representations},  (1999),
  preprint, Math.RT/9905045.

\bibitem{Ola-dg}
G.~\'Olafsson, \emph{Symmetric spaces of hermitian type}, Diff. Geom. and Appl.
  \textbf{1} (1991), 195--233.

\bibitem{oo-restric-1995}
G.~\'Olafsson and B.~\O{}rsted, \emph{Generalizations of the {Bargmann}
  transform}, Lie theory and its applications in physics. Proceedings of the
  international workshop, Clausthal, Germany, August 14-17, 1995. (H.-D.Doebner
  et~al, ed.), World Scientific, Singapore, 1996, pp.~3--14.

\bibitem{Opdam-bessel}
E.~Opdam, \emph{{Dunkl operators, Bessel functions and the discriminant of a
  finite Coxeter group}}, Compositio Math. \textbf{85} (1993), 333--373.

\bibitem{oz-weyl}
B.~\O{}rsted and G.~Zhang, \emph{Weyl quantization and tensor products of
  {Fock} and {Bergman} spaces}, Indiana Univ. Math. J. \textbf{43} (1994),
  551--582.

\bibitem{jp+gkz-pl3}
J.~Peetre and G.~Zhang, \emph{A weighted {Plancherel} formula
  \uppercase\expandafter{\romannumeral 3}. {The case of a hyperbolic matrix
  domain}}, Collect. Math. \textbf{43} (1992), 273--301.

\bibitem{Schmid}
W.~Schmid, \emph{Die {Randwerte} holomorpher {Funktionen} auf hermitesch
  symmetrischen {R\zae{}umen}}, Invent. Math \textbf{9} (1969), 61--80.

\bibitem{Torossian-res}
C.~Torossian, \emph{{Une application des op\'erateurs de Dunkl au th\'eor\`eme
  de restriction de Chevalley}}, C. R. Acad. Sci. Paris Sér. I Math.
  \textbf{318} (1994), no.~no. 10, 895--898.

\bibitem{UU}
A.~Unterberger and H.~Upmeier, \emph{The {Berezin} transform and invariant
  differential operators}, Comm. Math. Phys. \textbf{164} (1994), 563--597.

\bibitem{vanDiejen-tams}
J.~F. van Diejen, \emph{Properties of some family of hypergeometric orthogonal
  polynomials in several variables}, Trans. Amer. Math. Soc. \textbf{351}
  (1999), no.~1, 233--270.

\bibitem{Dijk-Pevzner}
G.~van Dijk and M.~Pevzner, \emph{Berezin kernels and tube domains}, J. Funct.
  Anal., to appear.

\bibitem{Wallach-sha}
N.~Wallach, \emph{Polynomial differential operators associated with {Hermitian}
  symmetric spaces}, Representation Theory of {Lie} Groups and {Lie} Algebras,
  Proceedings of Fuji-Kawaguchiko Conference 1990 (S.~Sano. T.~Kawazoe,
  T.~Oshima, ed.), World Scientific Singapore, 1992, pp.~76--95.

\bibitem{Yan-inv-eig}
Z.~Yan, \emph{Differential operators and function spaces}, Several complex
  variables in China, Contemp. Math., vol. 142, Amer. Math. Soc., 1993,
  pp.~121--142.

\bibitem{gkz-invtoe}
G.~Zhang, \emph{Tensor products of weighted {Bergman} spaces and invariant
  {Ha-plitz} operators}, Math. Scand. \textbf{71} (1992), 85--95.

\bibitem{gz-invdiff}
\bysame, \emph{Invariant differential operators on hermitian symmetric spaces
  and their eigenvalues}, Israel J. Math. \textbf{119} (2000), 157--185.

\bibitem{gkz-bere-rbsd}
\bysame, \emph{Berezin transform on real bounded symmetric domains}, Trans.
  Amer. Math. Soc. \textbf{353} (2001), 3769--3787.

\end{thebibliography}
\newcommand{\noopsort}[1]{} \newcommand{\printfirst}[2]{#1}
  \newcommand{\singleletter}[1]{#1} \newcommand{\switchargs}[2]{#2#1}
\providecommand{\bysame}{\leavevmode\hbox to3em{\hrulefill}\thinspace}

\end{document}